\documentclass{article}

\usepackage[letterpaper,asymmetric,left=1.25in,right=1.25in,top=1.25in,bottom=1.25in,
bindingoffset=0.0in]{geometry}

\usepackage{amsfonts,amsmath,amsthm,graphicx,color,epsfig}
\newtheorem{thm}{Theorem}[section]
\newtheorem*{Thm}{Theorem}
\newtheorem{prop}[thm]{Proposition}
\newtheorem{cor}[thm]{Corollary}
\newtheorem{lem}[thm]{Lemma}

\newtheorem{defn}[thm]{Definition}

\newcommand{\Endproof}{$\Box$ \vspace{.05in}}

\newcommand{\re}{{\rm Re}}
\newcommand{\size}{{\rm size}}
\newcommand{\vertline}{V}

\newcommand{\odd}{{\rm odd}}

\newcommand{\reg}{{\rm reg}}

\author{Roland K. W. Roeder}

\title{A degenerate Newton's Map in two complex variables: linking with currents}

\begin{document}

\maketitle

\begin{abstract} 
Little is known about the global structure of the basins of attraction of
Newton's method in two or more complex variables.  We make the first steps by
focusing on the specific Newton mapping to solve for the common roots of
$P(x,y) = x(1-x)$ and $Q(x,y) = y^2+Bxy-y$. 

There are invariant circles $S_0$ and $S_1$ within the lines $x=0$ and $x=1$
which are superattracting in the $x$-direction and hyperbolically repelling
within the vertical line.  We show that $S_0$ and $S_1$ have local
super-stable manifolds, which when pulled back under iterates of $N$ form
global super-stable spaces $W_0$ and $W_1$.  By blowing-up the points of
indeterminacy $p$ and $q$ of $N$ and all of their inverse images under $N$ we prove
that  $W_0$ and $W_1$ are real-analytic varieties.

We define linking between closed $1$-cycles in $W_i$ ($i=0,1$) and an
appropriate closed $2$ current providing a homomorphism
$lk:H_1(W_i,\mathbb{Z}) \rightarrow \mathbb{Q}$.  If $W_i$ intersects the
critical value locus of $N$, this homomorphism has dense image, proving that
$H_1(W_i,\mathbb{Z})$ is infinitely generated.  Using the Mayer-Vietoris exact
sequence and an algebraic trick, we show that the same is true for the closures
of the basins of the roots $\overline{W(r_i)}$.

\vspace{.05in}
{\footnotesize

}

\end{abstract}

Newton's method is one of the fundamental algorithms of mathematics, so it is
evidently important to understand its dynamics, particularly the structure of
the basins of attraction of the roots.  Even in one dimension, the topology of
these basins can be complicated and there has been a good deal of research on
this subject.  In higher dimensions, next to nothing is known about the
topology of the basins.  In this paper we make the first steps at understanding
their topology in two complex variables. 

We focus on a specific system: the Newton's Method used to solve for the common
roots of $P(x,y) = x(1-x)$ and $Q(x,y) = y^2+Bxy-y$.  While this is one
specific and relatively simple system, we believe that some of the techniques
developed in this paper can be used to study more general systems.

\vspace{.05in}

Dynamical systems $g:\mathbb{C}^n \rightarrow \mathbb{C}^n$ are often
classified in terms of:  (1) The number of inverse images of a generic point by
$g$, which is called the topological degree $d_t(g)$, and (2) Whether $g$ has
points of indeterminacy.

Mappings $g : \mathbb{P}^n \rightarrow \mathbb{P}^n$ with
$d_t(g)>1$, but without points of indeterminacy, have been studied
by Bonifant and Dabija \cite{Bonifant1}, Bonifant and Forn{\ae}ss \cite{Bonifant2}, Briend \cite{Br}, Briend and Duval \cite{BrD}, Dinh and Sibony \cite{DS2},
Fornaess and Sibony \cite{FS1,FS2,FS3}, Hubbard and Papadopol \cite{HP2},
Jonnson \cite{Jo}, and Ueda \cite{Ueda}.

Meanwhile, birational maps $g : \mathbb{P}^n \leadsto \mathbb{P}^n$
(rational maps with rational inverse) are examples of systems with points of
indeterminacy, but with $d_t(g)=1$.  The famous Henon mappings
$H:\mathbb{P}^2 \leadsto \mathbb{P}^2$ fall into this class.  Such systems
have been studied extensively by Bedford and Smillie
\cite{BS1,BS2,BS3,BS5,BS6,BS7,BS8}, Bedford, Lyubich and Smillie \cite{BS4},
Devaney and Nitecki \cite{DN}, Diller \cite{Dil}, Dinh and Sibony \cite{DS}, Dujardin \cite{Duj}, Favre and
Jonsson \cite{FJ}, Fornaess \cite{For}, Guedj \cite{Gu}, and Hubbard and
Oberste-Vorth \cite{HO1,HO2,HO3}.

Not nearly as much is known about mappings $g : \mathbb{P}^n \leadsto
\mathbb{P}^n$ with topological degree $d_t(g)> 1$ and with points of
indeterminacy.  The work of Russakovskii and Shiffman \cite{RS} considers a
measure that is obtained by choosing a generic point, taking the each of its
inverse images under $g^{\circ n}$ and giving them all equal weight in order to
obtain a probability measure $\mu_n$.  Under appropriate conditions on $g$ they
show that the measures $\mu_n$ converge to a measure $\mu$ that is independent
of the initial point.  In \cite{HP}, the authors present a proof by A. Douady
that $\mu$ does not charge points in the line at infinity, a result not
obtained in \cite{RS}.  In a recent paper, Guedj \cite{Gu2} shows that if the
topological degree $d_t(g)$ is sufficiently large, then $\mu$ does not charge
the points of indeterminacy of $g$ and does not charge any pluripolar set.  He
then uses these facts to establish ergodic properties of $\mu$.

\begin{figure}
\begin{center}
\scalebox{1.0}{
\begin{picture}(0,0)%
\epsfig{file=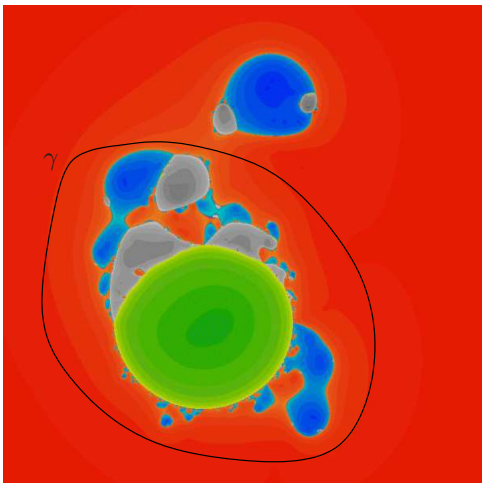}%
\end{picture}%
\setlength{\unitlength}{3947sp}%
\begin{picture}(3000,3000)(1201,-3361)
\put(1451,-1361){$\gamma$}%
\end{picture}%
}

\caption{Is the loop $\gamma$ trivial in  the homology of it's basin,
        within $\mathbb{C}^2$?}
\label{FIG_LOOP}
\end{center}
\end{figure}

Many of the papers considering mappings with both $d_t(g) > 1$ and points of
indeterminacy consider ergodic properties, invariant measures, and invariant
currents, focusing less on topological properties.  One paper that considers
some topological properties is \cite{HP}, by John Hubbard and Peter Papadopol,
who consider the dynamics of the Newton Map $N$ to solve for the zeros of two
quadratic equations $P$ and $Q$ in two complex variables.  The basins of
attraction for this system show interesting topology: for example, when drawing
intersections of a $1$ complex-dimensional slice with the basins of attraction
one often finds ``bubbles'' like the ones shown in Figure \ref{FIG_LOOP}.  It
is natural to ask if a loop, such as the one labeled $\gamma$ in the figure,
corresponds to a non-trivial loop in the homology of its basin of
attraction.  Clearly $[\gamma]$ is non-trivial in the basin intersected with
this slice, but it is much more difficult to determine if $[\gamma]$  is
non-trivial when considered within the entire 2 complex-dimensional basin,
which may reconnect in unusual ways outside of this slice.

Questions about the first homology of the basins are not answered by Hubbard
and Papadopol.   Using general principles they show that the basin of
attraction for each of the four roots is path connected and, by resolving points
of indeterminacy of $N$, they show that each basin is a Stein manifold.  But,
instead of addressing further questions about the topology of the individual
basins, they focus on creating a compactification with manageable topology on which $N$ is well-defined.
Most questions about the topology and
the detailed structure of these basins of attraction within their
compactification of $\mathbb{C}^2$ remain as mysteries.

In order to develop tools to answer some detailed questions about the topology of
basins of attraction for Newton maps $N:\mathbb{C}^2 \rightarrow \mathbb{C}^2$
to solve two quadratic equations $P(x,y)=0$ and $Q(x,y)=0$, we restrict our attention
to the degenerate case when the four roots lie on two parallel lines.
Normalizing, we study the Newton Map $N$ for $P(x,y) = x(1-x)$ and $Q(x,y)
= y^2+Bxy-y$.  In this case, the first component of $N(x,y)$ depends only on
$x$, while the second component depends on both $x$ and $y$.  Systems of this
form are commonly referred to as {\em skew products} in the literature
\cite{Ab,Hein,Jo,Ses,Sumi2,Sumi1} and they are often used as ``test cases''
when developing new techniques.  While we rely upon the fact that $N$ becomes a
skew product in this degenerate case, we hope that some of the techniques
developed here can eventually be adapted to non-degenerate cases.

\vspace{.05in}

Our approach is the following:

We compactify $\mathbb{C}^2$  obtaining a rational map $N:\mathbb{P}\times\mathbb{P}
\rightarrow \mathbb{P}\times\mathbb{P}$ with four points of indeterminacy at
$p=\left(\frac{1}{B},0 \right)$, $q =
\left(\frac{1}{2-B},\frac{1-B}{2-B}\right)$, $(\infty,\infty)$, and
$\left(\infty,\frac{B}{2}\right)$.  There are three invariant subspaces
$X_l:=\{(x,y) \mbox{ : } \re(x) < 1/2 \mbox{ and } x \neq \infty\},$ $X_{1/2}
:= \{(x,y) \mbox{ : } \re(x) = 1/2 \mbox{ or } x = \infty \}$ and $X_r:=
\{(x,y) \in  \mbox{ : } \re(x) > 1/2 \mbox{ and } x \neq \infty \}$.  The
common roots of $P$ and $Q$ are $r_1 = (0,0)$, $r_2=(0,1)$, $r_3 = (1,0)$, and
$r_4 = (1,1-B)$ with the basins of attraction $W(r_1)$ and $W(r_2)$ in $X_l$
and the basins of attraction $W(r_3)$ and $W(r_4)$ in $X_r$.  By
restricting to parameters $B \in \Omega = \{ B : |1-B| > 1\}$ we can assume
that both $p$ and $q$ are  in $X_l$.
 
We will prove that within $X_l$ and $X_r$ there are ``superstable
separatrices'' $W_0$ and $W_1$ consisting of the points that are attracted to
invariant circles within the lines $x=0$ and $x=1$ respectively.  By resolving
the points of indeterminacy of $N^k$ in $X_l$ we obtain a modified space
$X_l^\infty$ in which all iterates of $N$ are well-defined and in which $W_0$
is a real-analytic variety that provides a nice boundary between $W(r_1)$ and
$W(r_2)$.  Since $B \in \Omega$ there is no such problem in $X_r$: all
iterates of $N$ are already well defined on $X_r$ and $W_1$ is a real-analytic
variety in $X_r$.

In this paper we will study the topology of $W_0$ and $W_1$ in detail and we
will use a Mayer-Vietoris decomposition to relate their homology to the
homology of the basins of attraction of the four roots: $W(r_1),$ $W(r_2),$
$W(r_3),$ and $W(r_4)$ and the homology of $X_l^\infty$ and $X_r$.

The major emphasis of this paper is to show that loops in $W_0$ and $W_1$
that are generated by intersections of $W_0$ or $W_1$ with the critical value
locus $C$ are actually homologically non-trivial.  The
essential difficulty is to choose a notion of linking that is well defined
within the space $X_l^\infty$, which is very topologically complicated as a
result of the blow-ups.

We define linking between closed $1$-cycles in $W_i$ ($i=0,1$) and an
appropriate closed $2$ current providing a homomorphism
$lk:H_1(W_i,\mathbb{Z}) \rightarrow \mathbb{Q}$.  If $W_i$ intersects the
critical value locus of $N$, this homomorphism has dense image, proving that
$H_1(W_i,\mathbb{Z})$ is infinitely generated.  Using the Mayer-Vietoris exact
sequence and an algebraic trick, we show that the same is true for the closures
of the basins of the roots $\overline{W(r_i)}$.

Our work culminates to prove:

\begin{thm}\label{MAIN_THM}
Let $\overline{W(r_1)}$ and $\overline{W(r_2)}$ be the closures in $X_l^\infty$
of the  basins of attraction of $r_1=(0,0)$ and $r_2=(0,1)$ under iteration of
$N$ and let $\overline{W(r_3)}$ and $\overline{W(r_4)}$ be the closures in
$X_r$ of the basins of attraction of $r_3=(1,0)$ and $r_4=(1,1-B)$.
\begin{itemize} 

\item{$H_1\left(\overline{W(r_1)}\right)$ and $H_1\left(\overline{W(r_2)}\right)$ are infinitely generated for every $B \in \Omega$.}

\item{For $B \in \Omega$, if $W_1$ intersects the critical value parabola
$C(x,y)=0$, then both $H_1\left(\overline{W(r_3)}\right)$ and
$H_1\left(\overline{W(r_4)}\right)$ are infinitely generated,
otherwise $H_1\left(\overline{W(r_3)}\right)$
and $H_1\left(\overline{W(r_4)}\right)$ are trivial.}

\end{itemize}
\end{thm}

\vspace{.05in}

For $B \in \Omega_\reg$, the set of parameters for which
the separatrices are genuine manifolds, the basins of the four roots and
their closures in $X_l^\infty$ and $X_r$ have the some homotopy type.  Hence: 

\begin{cor}
For $B \in \Omega_\reg$, Theorem \ref{MAIN_THM} remains true when replacing
the closures of each of the basins with the basins themselves.
\end{cor}

\section{Basic properties of $N$}\label{BASICS}

In the first part of this section we summarize the basic results from \cite{HP}.

Given two vector spaces $V$ and $W$ of the same dimension and a differentiable mapping $F: V
\rightarrow W$, the associated Newton map $N_F: V \rightarrow V$ is given by
the formula
\begin{eqnarray}\label{NEWTON}
N_F({\bf x}) = {\bf x} - [DF({\bf x})]^{-1}(F({\bf x})).
\end{eqnarray}
\noindent

If $DF(r_i)$ is invertible for each root $r_i$ of $F$, then the roots of $F$
correspond to super attracting fixed points of $N_F$.  Conversely, every
fixed point of $N_F$ is a root of $F$.  Since each fixed point $r_i$ of $N_F$
is super-attracting, there is some neighborhood $U_i$ of $r_i$ for which each
initial guess ${\bf x_0} \in U_i$ will converge to $r_i$.  An explicit  lower
bound on the size of $U_i$ is given by Kantorovich's Theorem \cite{KANT}.

\begin{prop}{\rm \bf (Transformation rules)} \label{COORDCHANGE}
If $A: V \rightarrow V$ is affine, and invertible, and if $L:W \rightarrow W$ is linear and invertible,
then:
\begin{eqnarray}
N_{L\circ F\circ A} = A^{-1} \circ N_F \circ A.
\end{eqnarray}
\end{prop}
The proof is a careful use of the chain rule, see
\cite{HP}, Lemma 1.1.1.

\vspace{.05in}
\begin{prop}\label{DETERMINED}
Newton's Method to find the intersection of two quadratics depends only on the
intersection points and not on the choice of curves.
\end{prop}
\noindent For the proof, see Corollary 1.5.2 from \cite{HP}.

In this paper, we
normalize so that the roots are at ${0 \choose 0} $, ${1 \choose 0}$, ${0
\choose 1}$, and ${\alpha \choose \beta}$. If we let
$A=\frac{1-\alpha}{\beta}$ and $B = \frac{1-\beta}{\alpha}$, then $F {x \choose
y} = {x^2+Axy-x \choose y^2+Bxy-y} = {P(x,y) \choose Q(x,y)}$ has these roots and the corresponding
Newton Map is given by:
\begin{eqnarray} \label{NORMALIZATION}
N_F \left(\begin{array}{c} x \\ y \end{array} \right)   &=& \left(\begin{array}{c} x \\ y \end{array} \right) -
\left[ \begin{array}{cc} 2x+Ay-1 & Ax \\ By & 2y+Bx-1 \end{array} \right]^{-1} \left(\begin{array}{c} x^2+Axy-x \\ y^2+Bxy-y \end{array} \right)  \nonumber \\
&=& \frac{1}{\Delta}\left(\begin{array}{c} x(Bx^2+2xy+Ay^2-x-Ay) \\ y(Bx^2+2xy+Ay^2-Bx-y) \end{array} \right),
\end{eqnarray}
\noindent
where $\Delta = 2Bx^2+4xy +2Ay^2-(2+B)x-(2+A)y+1.$

\begin{prop}\label{HP_PROP1}
The critical value locus of $N_F$ is the union of the two parabolas
that go through the four roots of $F$.
\end{prop}

\begin{prop}\label{HP_PROP2}
The Newton Map has topological degree 4.
\end{prop}

See \cite{HP} for a proof of Propositions \ref{HP_PROP1} and \ref{HP_PROP2}.

\vspace{.05in}
It is a classical result that the dynamics of the Newton map
$N(z)$ to solve for the roots of any quadratic polynomial $p(z)$ is 
conjugate to the map $z \mapsto z^2$.  For the latter, the unit circle
$\mathbb{S}^1$ forms the boundary between the basin of attraction of $0$ and of
$\infty$.  If $\phi$ is the map conjugating $N(z)$ to $z \mapsto z^2$, then
$\phi^{-1}(\mathbb{S}^1)$ is the line in $\mathbb{C}$ that is equidistant from
the roots of $p$.  This line forms the boundary between the basin of the
two roots of $p(z)$ and the dynamics on this line (once you add a point at
infinity) are conjugate to angle doubling on the unit circle.

\begin{prop}{\rm ({\bf Invariant lines and invariant circles})}
\label{INVARIENT_LINES}
The lines joining the roots of $F$ are invariant under Newton Map $N_F$
and on these lines $N_F$ induces the dynamics of the one dimensional Newton
Map to find the roots of a quadratic polynomial.

Within each line is an invariant ``circle,'' corresponding to the points
of equal distance from the two roots in that line.
\end{prop}

\noindent  (See Proposition 1.5.3 in \cite{HP})

\noindent {\bf Proof:} Given any pair
of roots of $F$, there is an affine mapping taking them to ${ 0 \choose 0}$ and
${1 \choose 0}$ and a third root to ${0 \choose 1}$ The new system is also
normalized, but with the chosen pair of roots on
the $x$-axis.  Using Proposition \ref{COORDCHANGE}, if the $x$-axis is
invariant, then we will have shown that the line connecting the chosen pair of
roots is also invariant.  But this is easy to see because there is a factor of
$y$ in the second coordinate of Equation \ref{NORMALIZATION} for $N_F$.

The dynamics on the $x$-axis correspond to taking the first coordinate of $N_F$
in Equation \ref{NORMALIZATION} with $y=0$.  One finds $x \mapsto
\frac{x(Bx^2-x)}{2Bx^2-(2+B)x+1} = \frac{x^2}{2x-1}$.  This is the Newton's
Method to solve $x(1-x) = 0$.  Using the transformation rules from Proposition
\ref{COORDCHANGE} one can show that the same is true for any other invariant line.
\Endproof

\subsection{The degenerate case: $A=0$}\label{SUB_BASICS}
The
Newton map to find the common zeros of $P(x,y) = x(1-x)$ and $Q(x,y) = y^2+Bxy-y$ is:
\begin{eqnarray}\label{DEG_NORMALIZATION1}
N \left(\begin{array}{c} x \\ y \end{array} \right) = \frac{1}{\Delta}\left(\begin{array}{c} x(Bx^2+2xy-x) \\ y(Bx^2+2xy-Bx-y) \end{array} \right)
= \left(\begin{array}{c} \frac{x^2}{2x-1} \\ \frac{y(Bx^2+2xy-Bx-y)}{(2x-1)(Bx+2y-1)} \end{array} \right)
\end{eqnarray}

\noindent
with
\begin{eqnarray*}
\Delta = 2Bx^2+4xy-(2+B)x - 2y +1 = (2x-1)(Bx+2y-1).
\end{eqnarray*}
\noindent

\noindent
The fixed points of $N$ are the four common roots of $P$ and $Q$:
$r_1= (0,0), r_2 = (0,1), r_3 = (1,0)$, and $r_4 = (1,1-B)$.  

\begin{figure}[htp]
\begin{center}
\scalebox{.8}{
\begin{picture}(0,0)%
\epsfig{file=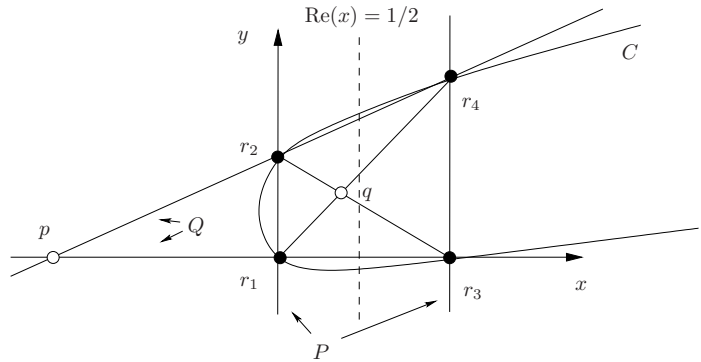}%
\end{picture}%
\setlength{\unitlength}{3947sp}%
\begin{picture}(5425,2856)(289,-2409)
\put(2612,291){${\rm Re}(x)=1/2$}
\put(2673,-2360){$P$}%
\put(2101,-736){$r_2$}%
\put(3858,-1863){$r_3$}%
\put(526,-1385){$p$}%
\put(5100, -2){$C$}%
\put(4734,-1822){$x$}%
\put(3061,-1089){$q$}%
\put(3843,-379){$r_4$}%
\put(2081,144){$y$}%
\put(1692,-1366){$Q$}%
\put(2101,-1786){$r_1$}%
\end{picture}%
}
\end{center}
\caption{The degenerate case $A=0$.}
\label{FIG_DEGCASE}
\end{figure}

The critical value locus is the union of the two parabolas going through the
four roots.  One of these coincides with $P(x,y) = x(1-x)=0$, while the other is
the non-degenerate parabola given by
\begin{eqnarray}\label{CRIT_VAL}
C(x,y) =y^2+Bxy+\frac{B^2}{4}x^2-\frac{B^2}{4}x-y=0.
\end{eqnarray}
We will often refer to the locus $C(x,y) = 0$ as {\em the critical value parabola} and denote it by $C$.
Figure \ref{FIG_DEGCASE} depicts the curves $P(x,y)=0$ and $Q(x,y)=0$,
the critical value parabola $C$, and the four roots, $r_1, r_2, r_3,$ and
$r_4$.

One can check directly from
Equation \ref{DEG_NORMALIZATION1} that $N$ has topological degree $4$, since every $x \neq 0,1$ has two
inverse images and the second component is an equation of degree two in $y$.

There are six invariant lines and, in this degenerate case,  these lines
have six points of intersection in $\mathbb{C}^2$.  Four of these intersections
correspond to the roots $r_1, r_2, r_3,$ and $r_4$, while the remaining two
correspond to points of indeterminacy. These are denoted $p$ and $q$ and are also
shown in Figure \ref{FIG_DEGCASE}.

The mapping governing the $x$ coordinate is $x \mapsto \frac{x^2}{2x-1}$, which
is itself the one variable Newton Map corresponding to the polynomial $x(x-1)$,
with Julia set consisting of the line $\re(x) = 1/2.$ This simple dynamics in $x$
is the main reason why the degenerate
Newton map is much easier to
understand than those considered in \cite{HP}: here all points in $\mathbb{C}^2$ with $\re(x) < 1/2$ are
super-attracted to the line $x=0$ and all points with $\re(x)>1/2$ are
super-attracted to the line $x=1$.  The vertical line at $x=m$ is mapped to the
line at $x=m^2/(2m-1)$ by the second coordinate of (\ref{DEG_NORMALIZATION1}),
which is in fact a rational map of degree 2, except at those values of $m$
where the numerator and the denominator in the second coordinate of
(\ref{DEG_NORMALIZATION1}) have a common factor.  This occurs exactly when $x=
1/B, x=1/(2-B),$ and $x=1/2.$ The first two correspond to the points of
indeterminacy $p$ and $q$.

There are three major invariant sets: $X_l:=\{(x,y) \mbox{ : } \re(x) < 1/2
\mbox{ and } x \neq \infty\},$ $X_{1/2} := \{(x,y) \mbox{ : } \re(x) = 1/2
\mbox{ or } x = \infty \}$ and $X_r:= \{(x,y) \in  \mbox{ : } \re(x) > 1/2
\mbox{ and } x \neq \infty \}$.  Figure \ref{FIG_DEGCASE} shows the case when
both points of indeterminacy $p$ and $q$ are in $X_l$.  The coordinates of $p$
and $q$ are $p = \left(\frac{1}{B},0\right)$ and
$q=\left(\frac{1}{2-B},\frac{1-B}{2-B}\right)$.  It is easy to check that $p$
and $q$ either are both in $X_l$, both in the separator $X_{1/2}$, or both in
$X_r$.  Let $\Omega = \{B \in \mathbb{C} \mbox{ : } |1-B|>1\}$ so that if $B
\in \Omega$ then both $p$ and $q$ are in $X_l$.  Using the transformation Rules
\ref{COORDCHANGE}, one sees that systems with this restriction still represent
every conjugacy class except for those corresponding to both $p,q \in X_{1/2}$.

\vspace{.05in}

\begin{figure}[htp]
\begin{center}
\scalebox{.8}{
\begin{picture}(0,0)%
\epsfig{file=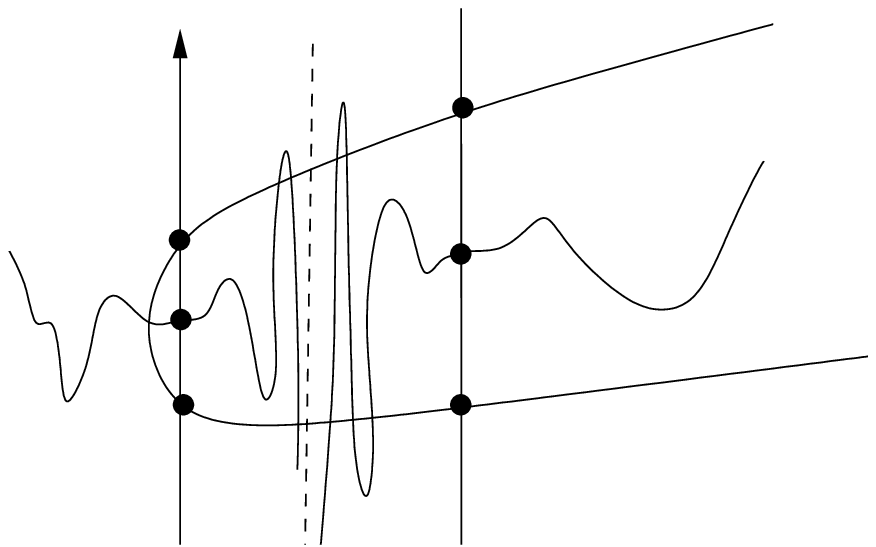}%
\end{picture}%
\setlength{\unitlength}{3947sp}%
\begin{picture}(4288,2695)(1426,-2248)
\put(2101,-1711){$r_1$}%
\put(3858,-1863){$r_3$}%
\put(5100, -2){$C$}%
\put(3843,-379){$r_4$}%
\put(2612,291){${\rm Re}(x)=1/2$}%
\put(2081,144){$y$}%
\put(2101,-736){$r_2$}%
\put(3826,-1036){$S_1$}%
\put(4276,-736){$W_1$}%
\put(1426,-1486){$W_0$}%
\put(2476,-1336){$S_0$}%
\end{picture}%

}
\end{center}
\caption{Superstable separatrices in the degenerate case, $A=0$.}
\label{FIG_DEGCASE_MANIFOLDS}
\end{figure}

Let $S_0$ and $S_1$ be the invariant circles in the fixed lines $x=0$ and $x=1$, respectively.
Because the lines $x=0$
and $x=1$ are super-attracting in the $x$-direction, $S_0$ and $S_1$
are super-attracting in the $x$-direction, as well. 
In Section \ref{MV} we will show that
these circles have local superstable manifolds $W^{loc}_0$ and $W^{loc}_1$. 
Pulling $W^{loc}_0$ and $W^{loc}_1$ back under the Newton map we generate
superstable spaces $W_0$ and $W_1$ that form the boundary between the basin
$W(r_1)$ and $W(r_2)$ and between the basin $W(r_3)$ and $W(r_4)$,
respectively.  Figure \ref{FIG_DEGCASE_MANIFOLDS} shows an illustration of
these separatrices. 

\begin{prop} \label{SYMMETRY} {\rm ({\bf Axis of symmetry})}
Let $\tau$ denote the vertical reflection about the line
$Bx+2y-1=0$, that is: $\tau(x,y) = (x,1-Bx-y)$.
Then, $\tau$ is a symmetry of $N$:
\begin{eqnarray*}
\tau \circ N = N \circ \tau.
\end{eqnarray*}
Furthermore, $N$ maps this axis of symmetry to the line $y= \infty$.
\end{prop}

\noindent
{\bf Proof:}
The map $\tau$ is affine and interchanges $r_1$ with $r_2$ and $r_3$ with
$r_4$.  Let $F {x \choose y} = {P(x,y) \choose Q(x,y)}$ so that $r_1, r_2,
r_3,$ and $r_4$ are the roots of $F$.  By Proposition \ref{DETERMINED},
the Newton map $N_{F \circ \tau}$ for finding the roots of
$F \circ \tau$ is the same as $N_F$, since they have the same roots.    By the
transformation rules of the Newton Map under affine coordinate changes,
$N_{F \circ \tau} = \tau^{-1} \circ N \circ \tau$.
Hence:
\begin{eqnarray*}
\tau \circ N = \tau \circ N_{F \circ \tau} =
\tau \circ \tau^{-1} \circ N \circ \tau = N \circ \tau
\end{eqnarray*}
\noindent
The axis $Bx+2y-1=0$ is mapped to the line $y=\infty$ due to the factor
$Bx+2y-1=0$ in the denominator the second component of $N$.  
\Endproof

\section{Computer exploration of $N$}\label{COMP_EXP}

In this section we show computer images of the basins of attraction for
the four common zeros of $P$ and $Q$ for the parameters $B=0.7857+1.1161i$, and
$B=-0.7902+1.7232i$.   All of the computer images displayed in this paper were made
using the wonderful program FractalAsm \cite{FA}.  

The separatrices $W_0$ and $W_1$ are clearly visible in these images, forming
the smooth boundary between $W(r_1)$ and $W(r_2)$ and between $W(r_3)$ and
$W(r_4)$, respectively.  The boundary between $W(r_1) \cup W(r_2)$ and $W(r_3)
\cup W(r_4)$, when visible, corresponds to points $(x,y)$ with $\re(x) = 1/2$.

\vspace{.05in}
\noindent
{\bf Case 1: $B=0.7857+1.1161i$}

\begin{figure}[htp]
\begin{center}
\scalebox{1.0}{
\begin{picture}(0,0)%
\epsfig{file=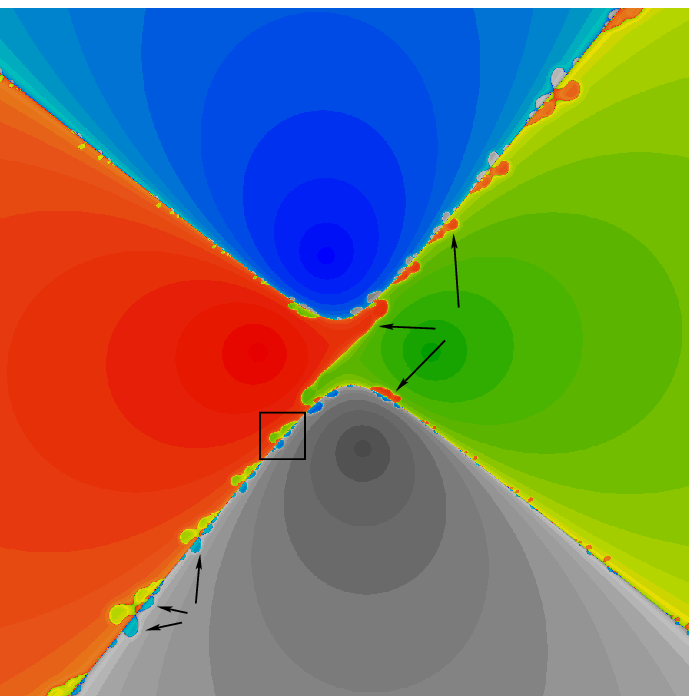}%
\end{picture}%
\setlength{\unitlength}{3947sp}%
\begin{picture}(3543,3304)(1201,-3665)
\put(1697,-2385){$W(r_1)$}%
\put(2565,-1022){$W(r_3)$}%
\put(2813,-3169){$W(r_4)$}%
\put(2102,-3345){$W_0 \cap C$}%
\put(3319,-1934){$W_0 \cap C$}%
\put(3745,-2398){$W(r_2)$}%
\end{picture}%

\begin{picture}(0,0)%
\epsfig{file=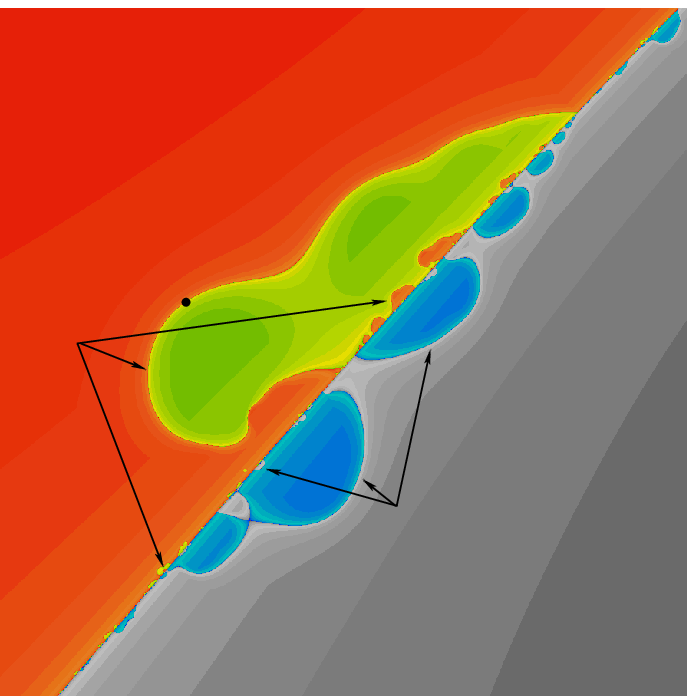}%
\end{picture}%
\setlength{\unitlength}{3947sp}%
\begin{picture}(3300,3300)(1201,-3661)
\put(3108,-2893){$W_1 \cap C$}%
\put(1967,-1705){$a$}%
\put(2556,-2532){$W(r_3)$}%
\put(2893,-1434){$W(r_2)$}%
\put(1614,-1063){$W(r_1)$}%
\put(2742,-3347){$W(r_4)$}%
\put(1247,-1897){$W_0 \cap C$}%
\end{picture}%

}
\end{center}
\caption{The critical value parabola $C$ for $B=0.7857+1.1161i$.
        The boundary between $W(r_1)$ and $W(r_2)$ is $W_0 \cap C$
        and the boundary between the $W(r_3)$ and $W(r_4)$ is $W_1 \cap C$.
        The image on the right is a zoomed in view of the boxed region in
        the image on the left.}
\label{FIG_PARAB1}
\end{figure}

The first kind of slice that we consider is given by the critical value
parabola $C$, which is
parameterized by a single complex variable, the offset from the axis of $C$.
Figure \ref{FIG_PARAB1} shows part of this slice on the left while the image
on the right shows a zoomed in view corresponding to the region enclosed in the
small rectangle from the image on the left.  The center of the symmetry $\tau$
is the center of the image on the left of Figure \ref{FIG_PARAB1}, but is
outside of the image on the right.

\begin{figure}[htp]
\begin{center}
\scalebox{.85}{
\begin{picture}(0,0)%
\epsfig{file=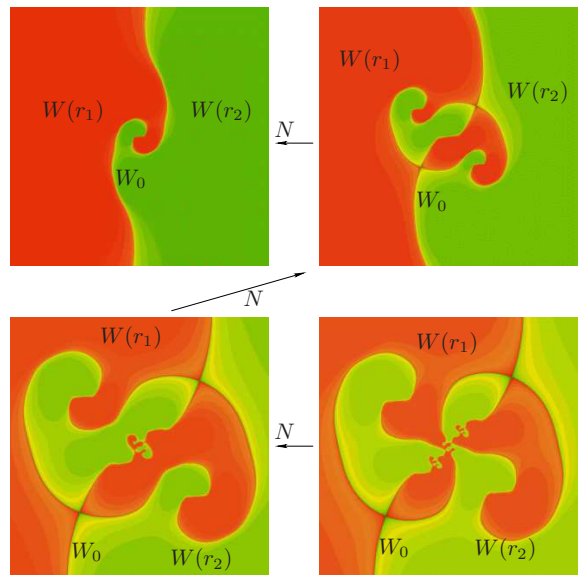}%
\end{picture}%
\setlength{\unitlength}{3947sp}%
\begin{picture}(4194,4194)(601,-4555)
\put(2555,-3555){$N$}%
\put(2555,-1314){$N$}%
\put(1936,-1171){$W(r_2)$}%
\put(4271,-1028){$W(r_2)$}%
\put(3031,-790){$W(r_1)$}%
\put(3603,-1838){$W_0$}%
\put(1268,-2840){$W(r_1)$}%
\put(4033,-4412){$W(r_2)$}%
\put(3603,-2887){$W(r_1)$}%
\put(1053,-4427){$W_0$}%
\put(1788,-4486){$W(r_2)$}%
\put(3323,-4401){$W_0$}%
\put(1379,-1673){$W_0$}%
\put(839,-1171){$W(r_1)$}%
\put(2328,-2566){$N$}%
\end{picture}%

}
\end{center}
\caption{Vertical line through point $a$ from Figure \ref{FIG_PARAB1} and
three inverse images of this line.  The boundary
between $W(r_1)$ and $W(r_2)$ is the intersection of $W_0$ with these
vertical lines.  Notice that there are many closed loops in $W_0$ within
these vertical lines.  The center of the symmetry $\tau$ is at the center of these
images.}
\label{FIG_PARAB1_SP1}
\end{figure}

Figure \ref{FIG_PARAB1_SP1} shows the vertical line $x=a$, where
$a$ is labeled in Figure \ref{FIG_PARAB1}, as well as the vertical lines
through three inverse images of $a$.  We have placed the center of the symmetry
$\tau$ at the center of these images so that reflection across this point
perfectly interchanges the basins.

Notice how the first inverse image of $x=a$ is divided into two regions that are in $W(r_1)$ and two
regions in $W(r_2)$.  
This is because we chose $a$ on the
superstable separatrix $W_0$.  
The lines at second and third inverse images of $x=a$ are divided
into three regions in $W(r_1)$ and in $W(r_2)$ and five regions in $W(r_1)$ and in $W(r_2)$, respectively.

\vspace{.05in}
\noindent
{\bf Case 2: $B=-0.7902+1.7232i$}

\begin{figure}[htp]
\begin{center}
\scalebox{1.0}{
\begin{picture}(0,0)%
\epsfig{file=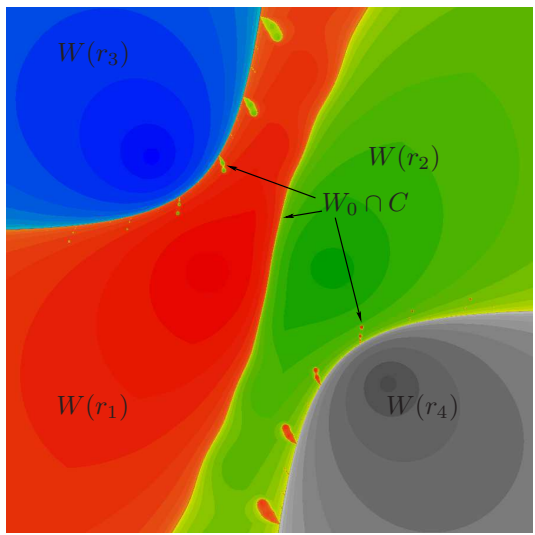}%
\end{picture}%
\setlength{\unitlength}{3947sp}%
\begin{picture}(3300,3300)(1201,-3661)
\put(3470,-1351){$W(r_2)$}%
\put(3181,-1640){$W_0 \cap C$}%
\put(3586,-2916){$W(r_4)$}%
\put(1513,-2917){$W(r_1)$}%
\put(1518,-699){$W(r_3)$}%
\end{picture}%

}
\end{center}
\caption{Critical value parabola $C$ for $B=-0.7902+1.7232i$.
The boundary between the $W(r_1)$ and $W(r_2)$ is $W_0 \cap C$.
We see no boundaries between $W(r_3)$ and $W(r_4)$, indicating
that $W_1$
might not intersect $C$.}
\label{FIG_PARAB2}
\end{figure}

Figure \ref{FIG_PARAB2} shows the intersections of the basins of attraction for
$W(r_1)$, $W(r_2)$, $W(r_3),$ and $W(r_4)$ with the critical value parabola
$C$.  Notice that there are clearly intersections of the superstable separatrix
$W_0$ with $C$ forming the visible boundary between $W(r_1)$ and $W(r_2)$.  However,
we see no boundaries between $W(r_3)$ and $W(r_4)$, indicating that $W_1$ might
not intersect $C$.  All of the further zoom-ins that we have done show no
evidence of intersections between $W_1$ and $C$.  We cannot prove that there
are values of $B$ for which $W_1 \cap C = \emptyset$, however it seems likely,
based on computer experiments.

\noindent

As for the previous value of $B$, the vertical lines above points of
intersection of $W_0$ with $C$ and the vertical lines mapped to them by $N$
contain many interestingly loops that are in $W_0$.

\section{Superstable separatrices $W_0$ and $W_1$.}\label{MV}

The invariant circle $S_0$ is the set of points in the line $x=0$ equidistant
from $r_1$ and $r_2$ and the invariant circle $S_1$ is the set of points in the
line $x=1$ equidistant from $r_3$ and $r_4$. 

\begin{prop}
The invariant circles $S_0$ and $S_1$ have multiplier $0$ in
the $x$-direction and they have multiplier $2$ within the vertical line in the direction normal
to the circle.
\end{prop}

\noindent
{\bf Proof:}
The vertical lines $x=0$ and $x=1$ are superattracting in the $x$-direction,
hence the circles $S_0$ and $S_1$ are as well.  Within these vertical lines,
$N$ is the Newton's method for the quadratic polynomial with roots $r_1$ and
$r_2$ (or $r_3$ and $r_4$), so the invariant circle is repelling in this line
with multiplier $2$.  \Endproof

\begin{prop}\label{LOCAL_STABLE_MFD}
The invariant circles $S_0$ and $S_1$ have local superstable manifolds
$W_0^{loc}$ and $W_1^{loc}$.
More specifically, there are neighborhoods $U_0, U_1 \subset \mathbb{C}$
of $x=0$ and $x=1$ and subsets $W_0^{loc} \subset
X_l$, $W_1^{loc} \subset X_r$
so that:

\begin{itemize}
\item $N(W_0^{loc}) \subset W_0^{loc}$ and $N(W_1^{loc}) \subset W_1^{loc}$

\item {$W_0^{loc}$ is the image of some $\Phi_0:U_0 \times S_0 \rightarrow
X_l$ which is analytic in the first coordinate and quasiconformal in the
second.}

\item {$W_1^{loc}$ is the image of some $\Phi_1:U_1 \times S_1 \rightarrow X_r$
which is analytic in the first coordinate and quasiconformal in the second.}

\end{itemize}
\end{prop}

\noindent
We use a technique due to John Hubbard and Sebastien Krief which allows us to
use the theory of holomorphic motions and the $\lambda$-Lemma of Ma\~ne,
Sad, and Sullivan \cite{MSS}, instead of the more standard graph transformation
approach.  A somewhat different stable manifold theorem for the invariant
circles in the non-degenerate case ($A \neq 0$) is also proved using this
technique in \cite{HP}.  While points in the manifolds obtained in our proof
are genuinely attracted to the circles $S_0$ and $S_1$, the situation in
\cite{HP} is much more complicated, with dense sets of points that are not
attracted to the invariant circles.

\vspace{.05in}
\noindent
{\bf Proof:}
To simplify computations we will make the change of variables $z(x) =
\frac{x}{x-1}$ and $w(y) = \frac{y}{y-1}$ which conjugates the first coordinate
of $N$ to $z \mapsto z^2$ and places the invariant circle $S_0$ at
$\{z=0, |w|=1\}$.  In the new coordinates $(z,w)$, the Newton map
becomes:
\begin{eqnarray} \label{NEW_COORDS}
N\left(\begin{array}{c} z \\ w \end{array} \right) =
\left(\begin{array}{c} z^2 \\
\frac{w^2+(Bw-Bw^2)z-w^2z^2}{1+(B-Bw)z+(Bw^2+B-1-2Bw)z^2} \end{array} \right)
\end{eqnarray}
and the critical value locus of $N$ in these coordinates is the image of
$C$ under the change of variables, which we denote by $C'$.  

Let
\begin{eqnarray*}
\Delta_{\epsilon,\delta} = \left\{ (z,w) \in X_l \mbox{ : } |z|<\epsilon \mbox{ and } 1 -\delta < |w| < 1+\delta \right\}
\end{eqnarray*}
\noindent
so that $\Delta_{\epsilon,\delta}$ is an open neighborhood of $S_0$.
The boundary of $\Delta_{\epsilon,\delta}$ consists of the vertical boundary
$\partial^V \Delta_{\epsilon,\delta} = \{|z| = \epsilon\}$ and the horizontal
boundary $\partial^H \Delta_{\epsilon,\delta} = \left\{
 |w| = 1 \pm \delta \right\}$.

We must choose $\epsilon$ and $\delta$ so that:
\begin{enumerate}
\item{$\Delta_{\epsilon,\delta}$ is disjoint from the critical value
        locus $C'$, and}

\item{ $N(\partial^H \Delta_{\epsilon,\delta})$ is entirely outside of
	$\Delta_{\epsilon,\delta}$ and $N(\partial^V \Delta_{\epsilon,\delta})$
	is entirely inside of $|z|<\epsilon.$}

\end{enumerate}

\begin{figure}[h!]
\begin{center}
\begin{picture}(0,0)%
\epsfig{file=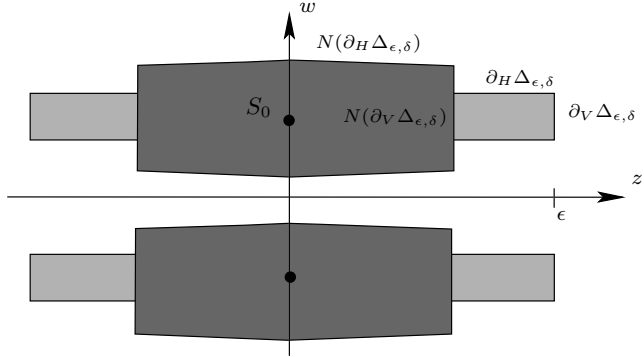}%
\end{picture}%
\setlength{\unitlength}{3947sp}%
\begingroup\makeatletter\ifx\SetFigFont\undefined%
\gdef\SetFigFont#1#2#3#4#5{%
  \reset@font\fontsize{#1}{#2pt}%
  \fontfamily{#3}\fontseries{#4}\fontshape{#5}%
  \selectfont}%
\fi\endgroup%
\begin{picture}(4116,2258)(151,-1617)
\put(2000,557){\makebox(0,0)[lb]{\smash{{\SetFigFont{8}{9.6}{\familydefault}{\mddefault}{\updefault}{\color[rgb]{0,0,0}$w$}%
}}}}
\put(3686,-105){\makebox(0,0)[lb]{\smash{{\SetFigFont{7}{8.4}{\familydefault}{\mddefault}{\updefault}{\color[rgb]{0,0,0}$\partial_V \Delta_{\epsilon,\delta}$}%
}}}}
\put(4083,-524){\makebox(0,0)[lb]{\smash{{\SetFigFont{8}{9.6}{\familydefault}{\mddefault}{\updefault}{\color[rgb]{0,0,0}$z$}%
}}}}
\put(2271,-133){\makebox(0,0)[lb]{\smash{{\SetFigFont{7}{8.4}{\familydefault}{\mddefault}{\updefault}{\color[rgb]{0,0,0}$N(\partial_V \Delta_{\epsilon,\delta})$}%
}}}}
\put(3166,111){\makebox(0,0)[lb]{\smash{{\SetFigFont{7}{8.4}{\familydefault}{\mddefault}{\updefault}{\color[rgb]{0,0,0}$\partial_H \Delta_{\epsilon,\delta}$}%
}}}}
\put(1661,-86){\makebox(0,0)[lb]{\smash{{\SetFigFont{9}{10.8}{\familydefault}{\mddefault}{\updefault}{\color[rgb]{0,0,0}$S_0$}%
}}}}
\put(2102,332){\makebox(0,0)[lb]{\smash{{\SetFigFont{7}{8.4}{\familydefault}{\mddefault}{\updefault}{\color[rgb]{0,0,0}$N(\partial_H \Delta_{\epsilon,\delta})$}%
}}}}
\put(3614,-753){\makebox(0,0)[lb]{\smash{{\SetFigFont{8}{9.6}{\familydefault}{\mddefault}{\updefault}{\color[rgb]{0,0,0}$\epsilon$}%
}}}}
\end{picture}%
\end{center}
\caption{Here $\Delta_{\epsilon,\delta}$ is shown
in light grey and  $N(\Delta_{\epsilon,\delta})$ in dark grey.}
\label{FIG_MFD}
\end{figure}

The critical value locus $C'$ intersects the vertical line $z=0$ transversely
at $w=0$ and $w=\infty$, so we can choose $\epsilon$ sufficiently small so that
$C'$ intersects $\mathbb{D}_\epsilon \times \mathbb{P}$ outside of
$\Delta_{\epsilon,\frac{1}{2}}$.  Now, we reduce $\epsilon$ and $\delta$ so
that the second condition holds.  Because the first coordinate of $N$ is $z
\mapsto z^2$, $N(\partial^V \Delta_{\epsilon,\delta})$ is automatically inside
of $|z| < \epsilon.$  In the line $z=0$, $N(z,w) = w^2$, so by continuity we
can choose $\epsilon$ and $\delta$ small enough that $N(\partial^H
\Delta_{\epsilon,\delta})$ is entirely outside of $\Delta_{\epsilon,\delta}$.

Let $\mathbb{D}_\epsilon$ be the open disc $|z|< \epsilon$ in $\mathbb{C}$ for this
$\epsilon$.  Conditions 1 and 2 on $\epsilon$ and $\delta$ were chosen so that
the following lemma is true:

\begin{lem}\label{DISC_PULLBACK}
Suppose that $D \subset \Delta_{\epsilon,\delta}$ is a complex disc that is
the graph of an analytic function $\eta: \mathbb{D}_\epsilon \rightarrow \mathbb{P}$.
Then $N^{-1}(D) \cap \Delta_{\epsilon,\delta}$ is the union of two disjoint
complex discs, each given as the graph of analytic functions
$\zeta_1, \zeta_2 : \mathbb{D}_\epsilon \rightarrow \mathbb{P}$.
\end{lem}

\noindent
{\bf Proof of Lemma \ref{DISC_PULLBACK}:}
The locus $N^{-1}(D) \cap \Delta_{\epsilon,\delta}$ satisfies the equation
$N(z,w) \in D$, which is equivalent to $N_2(z,w) = \eta(z^2)$, because $D$ is
the graph of $\eta$.  Since $D \subset \Delta_{\epsilon,\delta}$, $D$ is
disjoint from $C'$, so $\frac{\partial}{\partial w} N_2(z,w)$ is non-zero in a
neighborhood of $N^{-1}(D)$, and we can use the implicit function theorem to
solve for $w = \zeta_1(z)$ and $w=\zeta_2(z)$.  There are exactly two branches
because $N_2(z,w)$ is degree $2$ in $w$.

The graphs of $\zeta_1$ and $\zeta_2$ form the two complex discs
$N^{-1}(D) \cap \Delta_{\epsilon,\delta}$.
\Endproof Lemma \ref{DISC_PULLBACK}.

The line $w=1$ is invariant under $N$ and attracted to the point $(0,1) \in
S_0$.  Let $D_0 = \{(z,w) \mbox{ : } |z| < \epsilon, w=1 \}$.  We will form
$W_0^{loc}$ by taking inverse images of $D_0$.

Since $D_0 \subset \Delta_{\epsilon,\delta}$ satisfies the conditions of Lemma
\ref{DISC_PULLBACK}, letting $D_1 = N^{-1}(D_0) \cap \Delta_{\epsilon,\delta}$
we obtain two complex discs in $\Delta_{\epsilon,\delta}$ each of which is
given by the graph of an analytic function $\eta: \mathbb{D}_\epsilon \rightarrow
\mathbb{P}$ and each of which is mapped within $D_0$ by $N$.  These discs
intersect $S_0$ at $w=1$ and $-1$.

Because each of the discs in $D_1$ satisfies the hypotheses of Lemma
\ref{DISC_PULLBACK} we can repeat this process, letting $D_2 = N^{-1}(D_1) \cap
\Delta_{\epsilon,\delta}$, which this lemma guarantees is the union of four
disjoint discs in $\Delta_{\epsilon,\delta}$, each of which is the graph of
some analytic function $\eta: \mathbb{D}_\epsilon \rightarrow \mathbb{P}$.
These four discs intersect  $S_0$ at the fourth roots of $1$.  Repeating this
process, we obtain $D_n$ consisting of $2^n$ disjoint complex discs in
$\Delta_{\epsilon,\delta}$, each given by the graph of an analytic function
intersecting $S_0$ at the $2^n$-th roots of $1$.

Let
$D_\infty = \bigcup_{n=0}^\infty D_n$,
which consists of a union of disjoint complex discs through each of the dyadic
points ${\cal D}$ on $S_0$.  Each of these discs is the graph of an analytic function from
$\mathbb{D}_\epsilon$ to $\mathbb{P}$, and $D_\infty$ is forward
invariant to $S_0$ under $N$.

From a different perspective, $D_\infty$ prescribes a {\em
holomorphic motion}:
\begin{eqnarray*}
\phi: \mathbb{D}_\epsilon \times {\cal D} \rightarrow \mathbb{P}
\end{eqnarray*}
\noindent
where $\phi(z,\theta)$ is given by $\eta(z)$ where $\eta: \mathbb{D}_\epsilon
\rightarrow \mathbb{P}$ is the analytic function whose graph is the disc in
$D_\infty$ containing $\theta \in {\cal D} \subset S_0$.

By the $\lambda$-lemma of Ma\~ne-Sad-Sullivan \cite{MSS}, $\phi$ extends
continuously to a holomorphic motion on $S_0$, the closure of ${\cal D}$.
\begin{eqnarray*}
\phi: \mathbb{D}_\epsilon \times S_0 \rightarrow \mathbb{P}
\end{eqnarray*}

\noindent
We define $W_0^{loc}$ to be the image of $(z,w) \mapsto (z,\phi(z,w))$.
Clearly $N(W_0^{loc}) \subset W_0^{loc}$ and every point in
$W_0^{loc}$ is forward invariant to $S_0$.

The construction of of $W_1^{loc}$ is nearly identical and we omit it.
\Endproof Proposition \ref{LOCAL_STABLE_MFD}.

Because the local superstable manifolds $W_0^{loc}$ and $W_1^{loc}$ are
forward invariant under $N$, we can define global invariant
sets $W_0$ and $W_1$ by:
\begin{eqnarray*}
W_0 = \bigcup_{n=0}^\infty N^{-n}(W_0^{loc}), \qquad
W_1 = \bigcup_{n=0}^\infty N^{-n}(W_1^{loc}).
\end{eqnarray*}
\noindent

One might expect that $W_0$ and $W_1$ are manifolds, since the Inverse Function
Theorem gives that the pull-back $N^{-k}(W_i^{loc})$ (or $N^{-k}(W_i^{loc})$)
by $N$ is ``locally manifold'' at points where $N^{-(k-1)}(W_i^{loc})$ 
is disjoint from or transverse to the critical value
locus $C$.  However, we do expect that there will be some values of the
parameter $B$ for which there is a tangency between $N^{-k}(W_i^{loc})$
and $C$.  In fact, our computer images show that this must
be the case, because we see the topology of $W_0\cap C$ and of $W_1\cap C$
change as we change $B$.  For these parameter values $W_i$ will not
be a manifolds.  To make this distinction, we will call $W_0$ and $W_1$ {\em
separatrices} because they separate $W(r_1)$ from $W(r_2)$ and separate
$W(r_3)$ from $W(r_4)$.  

\vspace{.05in}
The following proposition requires that all iterates of $N$ be defined for all
points in a neighborhood of $W_0$ in $X_l$ and in a neighborhood of $W_1$ in
$X_r$.  This will require a modification $X_l^\infty$ of $X_l$ that is obtained
by blowing-up the points of indeterminacy $p$ and $q$ and all of their inverse
images under $N$.  We will prove Proposition \ref{GLOBAL_SEPARATRICIES},
temporarily thinking that we are working in $X_l$ and $X_r$, and then explain
why it is necessary to blow-up the points of indeterminacy.  
The entire construction of $X_l^\infty$ is given in the following section.

\begin{prop}\label{GLOBAL_SEPARATRICIES}
For every $B \in \Omega$, the separatrices $W_0$ and $W_1$ are real analytic subspaces of
$X_l^\infty$ and $X_r$, each defined as the zero set of a single
non-constant real-analytic equation in an neighborhood of $W_0$ and in
a neighborhood of $W_1$, respectively.
\end{prop}

\noindent
The proof is similar of that of B\"ottcher's Theorem in one
variable dynamics.
\vspace{.05in}

\noindent
{\bf Proof:}
We express $N$ in the variables $z = \frac{x}{x-1}$ and $w = \frac{y}{y-1}$
so that $S_0$ is given by $\{z=0, |w| = 1\}$.

We will show that
\begin{eqnarray*}
\phi(z,w) = \lim_{n \rightarrow \infty} (N^n_2(z,w))^{1/2^n}
\end{eqnarray*}
\noindent
converges to a non-constant analytic function on a
neighborhood of $W_0$.  Then, for every $(z,w) \in
W_0$, $|N^n_2(z,w)|$ converges to $1$ because $S_0 = \{|w|=1\}$, hence
$\omega(z,w) :=$ $\log|\phi(z,w)| =$ $\log|(N^n_2(z,w))^{1/2^n}|$ converges to $0$
on $W_0$ and to non-zero values away from $W_0$.

If $\phi$ converges, then $\phi$ and $\omega$ transform
nicely under the involution $\tau$.
For $|z|$ small and $|w|$ close to $1$, $\tau$ is close to $(z,w)
\mapsto (z,1/w)$.  Using this approximation, we have $\phi(\tau(z,w)) =
\lim_{n\rightarrow \infty}(N_2^n(\tau(z,w)))^{1/2^n} =
\lim(\tau(N^n(z,w))_2)^{1/2^n} \approx \lim (1/N_2^n(z,w))^{1/2^n} =
1/\phi(z,w).$  Consequently, $\omega(\tau(z,w)) = -\omega(z,w).$

\vspace{.05in}

We can write $\phi(z,w):=\lim_{n \rightarrow \infty} (N^n_2(z,w))^{1/2^n}$ as a
telescoping product:
\begin{eqnarray}\label{BIG_PRODUCT}
\phi(z,w) = N_2(z,w)^{1/2} \cdot
\frac{N^2_2(z,w)^{1/4}}{N_2(z,w)^{1/2}} \cdot
\frac{N^3_2(z,w)^{1/8}}{N^2_2(z,w)^{1/4}} \cdots
\end{eqnarray}

We now check that we can restrict the neighborhood of $W_0$ where
$\phi$ is defined so that we can use the binomial formula
\begin{eqnarray}\label{BINOM}
(1+u)^\alpha = \sum_{n=0}^\infty \frac{\alpha(\alpha-1)\cdots(\alpha-n+1)}
{n!} u^n, \mbox{    when }|u|<1
\end{eqnarray}
\noindent
to define the $\frac{1}{2^{n}}$-th root in the $n$-th term of this product.  We
do this first for points in a neighborhood of $W_0$ in $\overline{W(r_1)}$ and
a similar proof shows that the same works in a neighborhood of $W_0$ in
$\overline{W(r_2)}$.

\vspace{.05in}
In the coordinates $(z,w)$ the denominator of $N_2$ is of the form $1+r$ with
$r=(B-Bw)z+(Bw^2+B-1-2Bw)z^2$, so for $|(B-Bw)z+(Bw^2+B-1-2Bw)z^2| < 1$ the
second coordinate of $N$ can me written as:
\begin{eqnarray}\label{GOOD_FORM}
N_2(z,w) = w^2(1-Bz-z^2) + wz g(z,w)
\end{eqnarray}
\noindent
with $g(z,w)$ analytic.  We can write $N_2$ in form (\ref{GOOD_FORM}) in a
neighborhood of $z=0$ (hence a neighborhood of $S_0$) since
$|(B-Bw)z+(Bw^2+B-1-2Bw)z^2|$ vanishes when $z=0$.  From this point on, we
restrict our attention to this neighborhood of $z=0$.

\vspace{.05in}
\noindent
The general term $\frac{N^{n+1}_2(z,w)^{1/2^{n+1}}}{N^n_2(z,w)^{1/2^n}}$ in the product (\ref{BIG_PRODUCT}) is of the form
{\footnotesize
\begin{eqnarray*}
 \left(\frac{(N^n_2(z,w))^2(1-BN_1^n(z,w)-(N_1^n(z,w))^2)+N_2^n(z,w)N^n_1(z,w) \cdot  g\left(N^n_1(z,w),N^n_2(z,w)\right)}{(N^n_2(z,w))^2}\right)^{1/2^{n+1}} \cr
= \left(1-BN_1^n(z,w)-(N_1^n(z,w))^2+\frac{N_1^n(z,w)}{N^n_2(z,w)} \cdot
 g\left(N^n_1(z,w),N^n_2(z,w)\right) \right) ^{1/2^{n+1}}.
\end{eqnarray*}
}

We need to check that
we can restrict, if necessary, the neighborhood of definition for $\phi(z,w)$
so that
\begin{eqnarray}\label{LEQ}
\left|-BN_1^n(z,w)-(N_1^n(z,w))^2+\frac{N_1^n(z,w)}{N^n_2(z,w)} \cdot g\left(N_1^n(z,w),N^n_2(z,w)\right)\right|  \\ \leq
\left|BN_1^n(z,w)+(N_1^n(z,w))^2\right|+\left|\frac{N_1^n(z,w)}{N^n_2(z,w)} \cdot g\left(N_1^n(z,w),N^n_2(z,w)\right)\right|
\leq \frac{1}{2}.
\end{eqnarray}

The first term is not a problem because $N_1^n(z,w)=z^{2n}$ and we are
restricting to $|z|$ small.
Since we are in $\overline{W(r_1)}$ the only difficulty can
can occur if $N^n_2(z,w)$
goes to $0$ fast enough to make (\ref{LEQ}) large.  Detailed analysis of the behavior near $r_1$ resolves
this concern:

In \cite{HP}, the authors perform blow-ups at each of the four roots, and
observe that the Newton map $N$ induces rational functions of degree $2$ on
each of the exceptional divisors $E_{r_1}, E_{r_2}, E_{r_3},$ and $E_{r_4}$.
Let's  compute the rational function $s: E_{r_1} \rightarrow E_{r_1}$.
In the coordinate chart $m=\frac{z}{w}$, the extension to $E_{r_1}$
is obtained by:
\begin{eqnarray*}
s(m) = \lim_{w \rightarrow 0} \frac{m^2w^2(1+(B-Bw)mw+(Bw^2+B-1-2Bw)m^2w^2)}
{w^2+(Bw-Bw^2)mw-w^2m^2w^2} 
= \frac{m^2}{1+Bm},
\end{eqnarray*}
since $w=0$ on $E_{r_1}$.

The rational function $s(m)$ has $m=0$ as a superattracting fixed point, so
there is a neighborhood of $m=0 \in E_{r_1}$ within $\overline{W(r_1)}$ so that for
any point $(z,w)$ in this neighborhood, $\lim_{n\rightarrow \infty}
\left|\frac{N_1^n(z,w)}{N_2^n(z,w)}\right| = 0$.  Pulling back this neighborhood
under $N$ we find a neighborhood $V \subset \overline{W(r_1)}$ of the line $z=0$ in which
this limit holds.

Now we consider the case when $(z,w) \in \overline{W(r_2)}$.  The concern is that
$|w|$ may grow too fast for us to find a neighborhood of $S_0$ inequality
(\ref{LEQ}) is true.  Instead of analyzing the asymptotics of $g$, we
can re-write $N$ in the new coordinates $(z,s)$ with $s=1/w$ and a nearly identical
construction to that of $V$ gives the appropriate neighborhood $V'$. 

\vspace{.05in}

Restricting the points $(z,w) \in V \cup V'$, the $\frac{1}{2^{n+1}}$-th root in
the product (\ref{BIG_PRODUCT}) is well-defined.  We check that the product
converges on the neighborhood $\Lambda$ of $S_0$.  It is sufficient to show
that the corresponding series of logarithms converges.  The general term in
this series is:
{\small
\begin{eqnarray*}
\log \left| \left(1-BN_1^n(z,w)-(N_1^n(z,w))^2+\frac{N_1^n(z,w)}{N^n_2(z,w)} \cdot g\left(N_1^n(z,w),N^n_2(z,w)\right)
\right) ^{1/2^{n+1}} \right| \leq \frac{\log2}{2^{n+1}},
\end{eqnarray*}
}
\noindent
using Equation \ref{LEQ} and the triangle inequality so that
\begin{eqnarray*}
\left|1-
BN_1^n(z,w)-(N_1^n(z,w))^2+\frac{N_1^n(z,w)}{N^n_2(z,w)} \cdot g\left(N_1^n(z,w),N^n_2(z,w)\right)
\right| < 2.
\end{eqnarray*}

This sequence of logarithms converges because it is dominated by a geometric
series, and hence for the product (\ref{BIG_PRODUCT}) converges to the analytic
function on $\phi(z,w)$ on $\Lambda$.  This way $\omega(z,w) = \log|\phi(z,w)|$
is a real analytic function on $\Lambda$, and by the invariance properties of
$\phi$ on $\omega(z,w)$ is an analytic function on a neighborhood of $W_0$.

The proof that $W_1$ is the zero locus of a non-constant analytic
function is very similar.
\Endproof

{\em It is important to notice that in this proof we assumed that all iterates
of $N$ are defined at every point in $X_l$, forgetting temporarily the points
of indeterminacy $p$ and $q$ (and all of their inverse images in $X_l$.)} This
is a real problem because $W_0$ naturally goes through all of the points of
indeterminacy: Under a high enough iterate of $N$, the line
$x=\frac{1}{B(2-B)}$ is mapped by a ramified covering to a vertical line
arbitrarily close to the line $x=0$.  Since these lines intersect $W_0^{loc}$
in a topological circle, the line $x=\frac{1}{B(2-B)}$ intersects $W_0$ in a
(possibly more complicated) curve.  We will see that the resolution of the
indeterminacy at $p$ and $q$ replaces $p$ and $q$ with exceptional divisors
$E_p$ and $E_q$ that are mapped to $x=\frac{1}{B(2-B)}$ by isomorphisms.  So,
to make this proof correct, we will have to blow-up at $p$ and $q$, and, in fact,
at all inverse images of $p$ and $q$.

An alternative approach would be to study $W_0$ on $X_l - \cup_{n=0}^\infty
N^{-n}(\{p,q\})$, where we have already proven it is a real-analytic variety.
However, we want to consider the topology of $W_0$, without all of these points
removed, so we prefer to do the sequence of blow-ups.  \vspace{.05in}

\section{Resolution of points of indeterminacy}\label{COMP}

By restricting to parameters $B \in \Omega$, the points of
indeterminacy $p$ and $q$ are in $X_l$ and there are no points of indeterminacy
in $X_r$.  In this section we will describe how to resolve the indeterminacy in
$N$ at $p$ and $q$ and in higher iterates of $N$ at all of
the inverse images of $p$ and $q$ in $X_l$, obtaining a new space $X_l^\infty$ on which all
iterates of $N$ are defined at every point.  

Writing $W_0$ as a real-analytic variety is not the only motivation for the
construction of $X_l^\infty$.  We plan to study the detailed topology of the
basins of attraction and of the separatrices $W_0$ and $W_1$.  It is difficult to
decide what is a reasonable alternative to the statement of Theorem
\ref{MAIN_THM} without blowing up points.

\subsection{Construction of $X_l^\infty$ and $N_\infty:X_l^\infty \rightarrow 
X_l^\infty$.}\label{SECTION_BLOWUPS}

Most of the material in this section and in the following section closely 
follow the works of Hubbard and Papadopol \cite{HP} and Hubbard, 
Papadopol, and Veselov \cite{HPV}.

Substitution of the points $p$ and $q$ into $C(x,y)$ yields $\frac{1}{4}(B-1)$
and $\frac{B^2-7B+2}{4B-8}$, so for values of $B$ at which these expressions are
non-zero, neither $p$ nor $q$ is a critical value.

Let $S \subset \Omega$ be the subset of parameter 
space for which no inverse image of the point of indeterminacy $p$ or of 
point of indeterminacy $q$ is in the critical value locus $C$. 
We first describe the construction of $X_l^\infty$ for parameter 
values $B \in S$, and then explain the necessary modifications for special 
circumstance when $B \notin S$.

\begin{thm}\label{BAIRE} The set $S$ is generic in the sense of 
Baire's Theorem, i.e. uncountable and dense in $\Omega$.
\end{thm} 
Because of its computational nature, the proof of Theorem \ref{BAIRE} is in
Appendix \ref{AP_BAIRE}.

\vspace{.1in}
{\em \large Construction of $X_l^\infty$ when $B \in S$:}
\begin{prop} \label{PQ_BLOWUP}
Let $X_l^0$ be the space $X_l$ blown up at the points $p$ and $q$
and let $\pi_0:X_l^0 \rightarrow X_l$ be the corresponding projection.
\begin{itemize}
\item{The mapping $N$ extends analytically to a mapping $N_0:X_l^0 \rightarrow X_l$.}
\item{$N_0$ maps the exceptional divisors $E_p$ and $E_q$ to the line $x= \frac{1}{B(2-B)}$ by
	isomorphisms.}
\end{itemize}
\end{prop}

\noindent
{\bf Proof:}
The definition of a blow-up at a point is available
in Appendix \ref{APP_BLOW_UPS}.  Further details about blow-ups are available
in books on Algebraic Geometry, including \cite{GH}, and, in the context
of complex dynamics, in the papers
\cite{HPV,HP}.

We will work in the chart $(x,m) \mapsto (x,m(x-\frac{1}{B}),m) \in X_l \times
\mathbb{P}^1$.  
We write $N(x,y) =
(N_1(x,y),N_2(x,y))$ so that in the coordinates $(x,m)$ we have $N_1(x,m) =
\frac{1}{B(2-B)}$ and
{\small
\begin{eqnarray*}
N_2(x,m)  
= \frac{\frac{m}{B}(Bx^2+2xm(x-\frac{1}{B})-Bx-m(x-\frac{1}{B}))}{(2x-1)(\frac{2m}{B}+1)}
\end{eqnarray*}
}

\noindent 
When restricted to the exceptional divisor $E_p$ (set
$x=\frac{1}{B}$) the mapping becomes $m \mapsto \frac{m(1-B)}{(2-B)(2m+B)}.$ If
instead we had been working in the chart $(y,m') \mapsto
(m'y+\frac{1}{B},y,m')$, we would have obtained 
$m' \mapsto \frac{(1-B)}{(2-B)(2+m'B)}$.
This is consistent with the extension in terms of $m$ since one is obtained
from the other by the change of variables $m=\frac{1}{m'}$.  Both of the
expressions for $N$ restricted to $E_p$ are linear-fractional transformations,
hence $N$ maps $E_p$ to the line $y=\frac{1}{B(2-B)}$ by an isomorphism.

\vspace{.05in}
The exceptional divisor $E_q$ can be treated similarly.
\Endproof

We will denote the  vertical line $x=\frac{1}{B(2-B)}$ by $\vertline$,
since we use this line so frequently.  This is the vertical line that
is tangent to $C$ at its ``vertex''.

By construction, the extension  $N_0:X_l^0 \rightarrow X_l$ has no points of
indeterminacy.  However, since we need to iterate we must consider $N_0$ as
a rational map $N_0:X_l^0 \rightarrow X_l^0$.  Each of the inverse images of
$p$ and $q$ become points of indeterminacy of $N_0$ because we have blown up at
$p$ and $q$.
Because $B \in S$, neither $p$ nor $q$ are
critical values and each has four inverse images under $N_0$.   We can blow-up
at each of these eight points obtaining the space $X_l^1$ and the projection
$\pi_1:X_l^1 \rightarrow X_l^0$.  One can then extend $N_0$ to the exceptional
divisors, obtaining $N_1:X_l^1 \rightarrow X_l^0$.

To make iterates $N^{\circ k}$ of $N$ well-defined for all $k$ we must repeat 
this process for the $k$-th inverse images, obtaining successive blow-ups 
$\pi_k: X_l^k \rightarrow X_l^{k-1}$ for every $k$.  The following proposition 
describes the extension of $N$ to these spaces:

\begin{prop} \label{REPEATED_BLOWUP}
Denote by $X_l^k$ the space $X_l^{k-1}$ blown up at each of these $2 \cdot 
4^k$ $k$-th inverse images of $p$ and $q$.
\begin{itemize}
\item{The mapping $N_{k-1}$ extends analytically to a mapping 
	$N_k:X_l^k \rightarrow X_l^{k-1}$.}
\item{Suppose that $z$ is one of the $k$-th inverse images of $p$ or $q$ and
denote the exceptional divisor over $z$ by $E_z$.  Then, $N_k$ maps $E_z$ to
$E_{N(z)}$ by isomorphism.}
\end{itemize}
\end{prop}

\noindent
{\bf Proof:}
As in Proposition \ref{PQ_BLOWUP}, denote the first and second components of $N$ by
$N_1(x,y)$ and $N_2(x,y)$.  Then, in the coordinates
$(x,m) \mapsto (x,mx,m)$ in a neighborhood of $E_z$ the mapping is given by: $m
\mapsto \frac{\partial_x N_1 + \partial_y N_1 m} {\partial_x N_2 + \partial_y
N_2 m}$.
By the assumption that $B \in S$, $DN$ is non-singular at $z$ and this gives an isomorphism from $E_z$
to $E_{N(z)}$.
\Endproof

\vspace{.05in}
Hence, by repeated blow-ups  we obtain a sequence of spaces and projections:
\begin{eqnarray}\label{INVSEQ}
X_l \xleftarrow{\pi_0} X_l^0 \xleftarrow{\pi_1} X_l^1 \xleftarrow{\pi_2} X_l^2 \xleftarrow{\pi_3} X_l^3 \xleftarrow{\pi_4} X_l^4 \xleftarrow{\pi_5} X_l^5 \xleftarrow{\pi_6} \cdots
\end{eqnarray}
The extensions of the Newton map $N$ to these spaces that we calculated in 
Propositions \ref{PQ_BLOWUP} and \ref{REPEATED_BLOWUP} give another 
sequence of spaces and mappings:
\begin{eqnarray}\label{MAPSEQ}
X_l \xleftarrow{N_0} X_l^0 \xleftarrow{N_1} X_l^1 \xleftarrow{N_2} X_l^2 \xleftarrow{N_3} X_l^3 \xleftarrow{N_4} X_l^4 \xleftarrow{N_5} X_l^5 \xleftarrow{N_6} \cdots
\end{eqnarray}

However, we do not have a single space $X_l^\infty$, nor a single mapping 
$N_\infty$ from this space to itself.   There is a standard 
procedure using {\em Inverse Limits} to create such a space and mapping 
from a sequence of spaces (\ref{INVSEQ}) and the sequence of mappings like 
(\ref{MAPSEQ}). That is, we will let $X_l^\infty$ be the inverse limit of 
the blown up spaces and projections in Sequence \ref{INVSEQ} and then use 
the sequence of extensions of the Newton maps \ref{MAPSEQ} to define a 
mapping $N_\infty : X_l^\infty \rightarrow X_l^\infty$ which naturally 
corresponds to an extension of $N$.

\begin{defn}
An {\bf Inverse system}, denoted $(M_i, \sigma_i)$, is a
family of objects $M_i$ in a category $C$ indexed by the natural numbers and
for every $i$ a morphism $\sigma_{i}:M_i \rightarrow M_{i-1}$.

\vspace{.05in}
\noindent
The {\bf Inverse Limit} of an inverse system $(M_i, \sigma_i)$, denoted by 
$\varprojlim (M_i,\sigma_i)$, is the object $X$ in $C$ together with 
morphisms $\alpha_i: X \rightarrow M_i$ satisfying $\alpha_{i-1} = 
\sigma_i \circ \alpha_i$ for each $i$ that is determined uniquely by the following universal 
property:

For any other pair $Y,\beta_i:Y\rightarrow M_i$ such that $\beta_{i-1}=\sigma_i
\circ \beta_i$, we have a unique morphism $u:Y \rightarrow X$ so that for each
$i$ we have $\beta_i = \alpha_i \circ u$.  \end{defn}

For our uses, the category will always be analytic spaces and the morphisms
holomorphic maps.  While not needed here, inverse systems and
inverse limits can be defined more generally, for objects $M_i$ indexed by a filtering
partially ordered set $I$.  The following proposition gives a construction of
$\varprojlim (M_i,\sigma_i)$ as a subset of the product space $\Pi_i M_i$.

\begin{prop}\label{INVLIM2}
Given an inverse system $(M_i,\sigma_{i})$ indexed by $\mathbb{N}$ (i.e. 
$\sigma_i: M_i \rightarrow M_{i-1}$), we can construct the inverse limit as 
follows:
\begin{eqnarray*}
\varprojlim (M_i,\sigma_i) = \{(m_0,m_1,m_2,m_3,\cdots) | m_i \in M_i \mbox{ and }
\sigma_{i}(m_i) = m_{i-1}\}.
\end{eqnarray*}
\end{prop}

\vspace{.05in}

We define $X_l^\infty = \varprojlim (X_l^k,\pi_k)$.  Using Proposition 
\ref{INVLIM2} we can state more concretely that
\begin{eqnarray*}
X_l^\infty = \{(x_0,x_1,x_2,x_3,\cdots) | x_i \in X_l^i \mbox{ and }\
\pi_i(x_i) = x_{i-1}\}. 
\end{eqnarray*}
In this description of $X_l^\infty$, the mappings $N_k:X_l^k \rightarrow X_l^{k-1}$ induce a mapping $N_\infty : X_l^\infty
\rightarrow X_l^\infty$ given by
$N_\infty((x_0,x_1,x_2,x_3,\cdots)) = (N_1(x_1),N_2(x_2),N_3(x_3),\cdots).$

\vspace{.1in}
{\em \large Construction of $X_l^\infty$ when $B \notin S$:\\} 
For parameter
values $B \notin S$, the blow-ups done at $p$ and $q$ in Proposition
\ref{PQ_BLOWUP} are exactly the same, since $N$ extends to
these blow-ups for any value of $B$.  (It is worth noticing that there is
actually a critical point of $N$ on both $E_p$ and on $E_q$ since each is
mapped to the line $V$, which contains a point of $C$.)  Special care
needs to be taken when a $k$-th inverse image of $p$ and of $q$ is a critical
point of $N$.  We describe the process here, although we leave
some of the details for the appendix.

The goal is to produce a space $X_l^k$ and a projection $\pi_k:X_l^k
\rightarrow X_l^{k-1}$ in such a way that $N$ extends to a map (without
singularities) $N_k:X_l^k \rightarrow X_l^{k-1}$.  If we can create the spaces
$X_l^k$ and extensions $N_k$ at every ``level'' $k$, we can use the
same inverse sequence process to make $X_l^\infty$ and $N_\infty:X_l^\infty \rightarrow
X_l^\infty$.

\begin{figure}
\begin{center}
\scalebox{.8}{
\begin{picture}(0,0)%
\epsfig{file=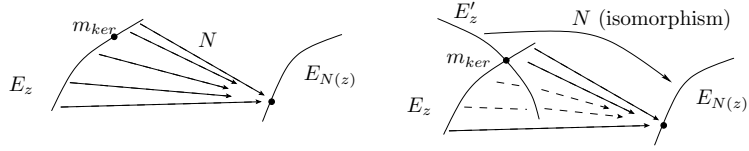}%
\end{picture}%
\setlength{\unitlength}{3947sp}%
\begingroup\makeatletter\ifx\SetFigFont\undefined%
\gdef\SetFigFont#1#2#3#4#5{%
  \reset@font\fontsize{#1}{#2pt}%
  \fontfamily{#3}\fontseries{#4}\fontshape{#5}%
  \selectfont}%
\fi\endgroup%
\begin{picture}(6169,1139)(400,-493)
\put(903,419){\makebox(0,0)[lb]{\smash{{\SetFigFont{10}{12.0}{\familydefault}{\mddefault}{\updefault}{\color[rgb]{0,0,0}$m_{ker}$}%
}}}}
\put(5794,-157){\makebox(0,0)[lb]{\smash{{\SetFigFont{10}{12.0}{\familydefault}{\mddefault}{\updefault}{\color[rgb]{0,0,0}$E_{N(z)}$}%
}}}}
\put(3892,521){\makebox(0,0)[lb]{\smash{{\SetFigFont{10}{12.0}{\familydefault}{\mddefault}{\updefault}{\color[rgb]{0,0,0}$E'_z$}%
}}}}
\put(4822,467){\makebox(0,0)[lb]{\smash{{\SetFigFont{10}{12.0}{\familydefault}{\mddefault}{\updefault}{\color[rgb]{0,0,0}$N$ (isomorphism)}%
}}}}
\put(3526,-217){\makebox(0,0)[lb]{\smash{{\SetFigFont{10}{12.0}{\familydefault}{\mddefault}{\updefault}{\color[rgb]{0,0,0}$E_z$}%
}}}}
\put(3838,161){\makebox(0,0)[lb]{\smash{{\SetFigFont{10}{12.0}{\familydefault}{\mddefault}{\updefault}{\color[rgb]{0,0,0}$m_{ker}$}%
}}}}
\put(2716, 29){\makebox(0,0)[lb]{\smash{{\SetFigFont{10}{12.0}{\familydefault}{\mddefault}{\updefault}{\color[rgb]{0,0,0}$E_{N(z)}$}%
}}}}
\put(1888,293){\makebox(0,0)[lb]{\smash{{\SetFigFont{10}{12.0}{\familydefault}{\mddefault}{\updefault}{\color[rgb]{0,0,0}$N$}%
}}}}
\put(400,-56){\makebox(0,0)[lb]{\smash{{\SetFigFont{10}{12.0}{\familydefault}{\mddefault}{\updefault}{\color[rgb]{0,0,0}$E_z$}%
}}}}
\end{picture}%

}
\end{center}
\caption{Blowing up a point on an exceptional divisor.}
\label{FIG_SECOND_BLOW_UP}
\end{figure}

Suppose for the moment that $z$ is a $k$-th inverse image of $p$ and that none
of the $n$-th inverse images of $p$ for $n<k$ were in the critical locus
$N^{-1}(C)$.  In this case, there is a single exceptional divisor in
$X_l^{k-1}$ above $N(z)$.  Because $z$ is critical, the extension of $N$ to
$E_z$ will map all of $E_z$ (except for one point) to a single point in
$E_{N(z)}$.  However, at
the slope $m_{ker} \in E_z$ corresponding to the kernel of $DN$, the extension to
$E_z$ has another point of indeterminacy!  Consequently, one has to blow-up
this point on $E_z$, obtaining a second exceptional divisor $E'_z$ above
$m_{ker}$.  Figure \ref{FIG_SECOND_BLOW_UP} shows this situation.  
An easy check using Taylor series for $N$ shows that
$N$ extends to $E'_z$ giving isomorphism from $E'_z$ to $E_{N(z)}$.

These two blow-ups above $z$ are sufficient to extend $N$.  

The two exceptional divisors above $z$ result in a further complication at
every point $w$ that is mapped to $z$.  Suppose that we have blown up at $w$.
The extension of $N$ to $E_w$ has a point of indeterminacy at the point that is
mapped to $m_{ker}\in E_z$.  Because of this, one has to blow-up a second time
above $w$ to resolve this point of indeterminacy.  In fact, at every repeated
inverse image of $z$ one will have to blow-up at least twice to resolve $N$.

There are further problems if an inverse image of $z$ is again critical.
At such a point, one will have to do even more blow-ups to resolve $N$.
A detailed description of this process becomes rather tedious, and we
will stop here.

\subsection{The mappings from $E_z$ to $\vertline$}
We saw in the previous section that $N$ maps each exceptional divisor that 
was newly created in $X_l^k$ to one of the exceptional divisors newly 
created in $X_l^{k-1}$ by either an isomorphism, or a constant map.
Since $N$ maps each $E_p$ and $E_q$ isomorphically to the line 
$\vertline$, the composition $N^{\circ k+1}$ maps each of the newly 
created exceptional divisors $E_z$ in $X_l^k$ to $\vertline$ either
by an isomorphism, or a constant map.  In summary:

\begin{prop}\label{EXCEPIONAL_COVERS}
Let $E_z$ be one of the exceptional divisors newly created in $X_l^k$
and let $\vertline$ be the line $x=1/(B(2-B))$.
Then $N^{\circ k+1}$ maps $E_z$ to $\vertline$ by an isomorphism, or a
constant map.
\end{prop}

\subsection{Homology of $X_r$ and of $X_l^\infty$}
Our goal in this paper is to relate the homology of the basins of attraction
for the four roots of $F$ to the homology of the spaces $X_r$ and $X_l^\infty$
and to the homology of the separatrices $W_0$ and $W_1$.  In this section, we
will compute the homology of $X_r$ and $X_l^\infty$.

Given a set $T$, we will denote by $\mathbb{Z}^{(T)}$ the submodule of the
product $\mathbb{Z}^{T}$ for which each element has at most finitely many non-zero
components.

We often find it necessary to encode information about the generators 
homology within our notation.  For example, the
module $\mathbb{Z}^{\{[K]\}}$ means the module  $\mathbb{Z}$ that is generated
by the fundamental class of $[K]$.

\begin{prop} 
We have: $H_0(X_r) = \mathbb{Z}$, $H_2(X_r) = \mathbb{Z}^{\{[\mathbb{P}^1]\}}$,
and $H_i(X_r) = 0$, for $i \neq 0$ or $2$.
\end{prop}

Unfortunately homology does not behave nicely under inverse limits.  So,
instead of directly using the fact that $X_l^\infty$ is an inverse limit to
compute its homology, we will write $X_l^\infty$ is a union of open subsets
$U_0 \subset U_1 \subset U_2 \subset \cdots$ in such a way that $H_2(U_i) =
\mathbb{Z}^{(L_i \cup \{[\vertline ]\})}$, where $L_i$ is the set of fundamental
classes of exceptional divisors contained in $U_i$ and $[V]$ is the fundamental
class of the vertical line  $\vertline$ given by $x=\frac{1}{B(2-B)}$.

Recall that the projection $\pi: X_l^\infty \rightarrow X_l$ is continuous, we
will create an exhaustion of $X_l^\infty$ by open sets $U_0 \subset U_1 \subset
U_2 \subset \cdots$ that are inverse images of open subsets forming an
exhaustion of  $X_l$.
Let $V_k = X_l - \bigcup_{n=k}^{\infty}\{N^{-n}(p),N^{-n}(q)\}$.  Clearly
$V_k$ is an open subset of $X_l$, so we will let $U_k = \pi^{-1}(V_k)$.
It is also clear that $U_1 \subset U_2 \subset U_3 \subset \cdots $
and that $\bigcup_{k=1}^\infty U_k = X_l^\infty$.

\begin{lem} For each $k$, $H_2(U_k) \cong H_2(X_l^k)$
\end{lem}

\noindent
{\bf Proof:}
Notice that $U_k$ canonically isomorphic to $X_l^k - \bigcup_{n=k}^{\infty}\{N^{-n}(p),N^{-n}(q)\}$.  Removing a discrete set of points from a $4$ (real)
dimensional manifold does not affect the second homology.  Hence,
$H_2(U_k) \cong H_2(X_l^k)$.
\Endproof

\begin{lem} $H_2(X_l^k) \cong \mathbb{Z}^{(L_k \cup \{[\vertline]\})}$, where
$L_k$ is the set of fundamental classes of exceptional divisors in $X_l^k$.
\end{lem}

\begin{prop}\label{HOMOLOGY} $H_2(X_l^\infty) \cong \mathbb{Z}^{(L \cup \{[\vertline]\})}$, where
$L$ is the set of fundamental classes of exceptional divisors in $X_l^\infty$ and $[\vertline]$
is the fundamental class of the vertical line $\vertline$.
\end{prop}

\noindent
{\bf Proof:}
Since $X_l^\infty = \bigcup_{k=1}^\infty U_k$ and $H_2(U_k) \cong H_2(X_l^k)
\cong \mathbb{Z}^{(L \cup \{[\vertline]\})}$, we 
have that $H_2(X_l^\infty) \cong \varinjlim \left( \mathbb{Z}^{(L_k \cup \{[\vertline]\})}\right)$, which is 
$\mathbb{Z}^{(L \cup \{[\vertline]\})}$.
\Endproof

In the generic case $B \in S$ for which none of the inverse images of $p$ or $q$ under
$N$ are in the critical value parabola $C$ we can describe $H_2(X_l^\infty)$
somewhat more explicitly:
\begin{prop}\label{HOMOLOGY2}
For $B \in S$,
$H_2(X_l^\infty) = \mathbb{Z}^{\{[\vertline]\}} \oplus \left(\bigoplus_{N^k(x) = p} \mathbb{Z}^{\{[E_{x}]\}}\right) \oplus \left( \bigoplus_{N^k(x) = q} \mathbb{Z}^{\{[E_{x}]\}}\right)$.
\end{prop}

\noindent
{\bf Proof:}
This is merely a restatement of Proposition \ref{HOMOLOGY} using that when
$B \in S$, only a single blow-up is necessary at each $k$-th inverse
image of $p$ and of $q$ for every $k$.
\Endproof

We will need the following proposition about the intersection of classes
in $H_2(X_l^\infty)$:

\begin{prop}\label{SELF_INTERSECTION}
Let $[\vertline]$ and $[E_z]$ be the fundamental classes
of a vertical line $\vertline$ and an exceptional divisor $E_z$ in
$H_2(X_l^\infty)$. Then
$[\vertline] \cdot [\vertline] = 0$ and
$[E_z] \cdot [E_z] \leq -1$.
\end{prop}

\noindent
{\bf Proof:}
We have chosen the vertical line $\vertline$ so that points on it are never
blown up, hence within $X_l^\infty$ it has self-intersection number $0$, just
as it did in $X_l$.  

If no points on the exceptional divisor $E_z$ have been blown up,
then it is a classical result that $[E_z]\cdot[E_z] = -1$.  Otherwise,
if points in $E_z$ have been blown up, it is a classical result that 
each blow-up reduces  $[E_z]\cdot[E_z]$ by $1$, hence $[E_z]\cdot[E_z]
\leq -1$.  (See \cite{GH}.)
\Endproof

\subsection{Why we work in $X_l$}

Suppose for a moment that we did this sequence of blow-ups in $X =
\mathbb{P}^1\times \mathbb{P}^1$, instead of just within the invariant subspace
$X_l$.  The repeated inverse images of the points of indeterminacy must
accumulate because $X$ is compact.  The topology of the inverse limit
$X^\infty$ becomes very complicated near these points of accumulation.  Hubbard
and Papadopol develop elaborate techniques including {\em Farey Blow-ups} and
{\em Real-oriented Blow-ups}  to ``tame'' the topology at these points of
accumulation.  We avoid this problem caused by accumulation  by working in the
space $X_l$ since the inverse images of $p$ and $q$ go to the ``end'' of $X_l =
\{\re(x) < 1/2\} \times \mathbb{P}^1$ instead of accumulating to a finite
point.

\section{Mayer-Vietoris sequences}
We will study the topology
of $W_0$ and $W_1$ in detail in order to prove Theorem \ref{MAIN_THM}.  The following Mayer-Vietoris calculations
will allow us to relate their homology to that of the basins of
attraction for the four roots $r_1,$ $r_2$, $r_3$, and $r_4$.

Let $\overline{W(r_1)}$ and $\overline{W(r_2)}$ be the closures of $W(r_1)$ and
$W(r_2)$ in $X_l^\infty$ and let $\overline{W(r_3)}$ and $\overline{W(r_4)}$ be
the closures of $W(r_3)$ and $W(r_4)$ in $X_r$.  Since $W_0$ and $W_1$ are
real-analytic varieties in $X_l^\infty$ and $X_r$, respectively, there are
neighborhoods in $X_l^\infty$ and $X_r$ of $W_0$ and $W_1$ that deformation
retract onto $W_0$ and $W_1$.  Hence, we can consider the Mayer-Vietoris exact sequence
(see \cite{BREDON, HATCHER}) for the decompositions $\overline{W(r_1)} \cup
\overline{W(r_2)} = X_l^\infty, \mbox{    }$ $\overline{W(r_1)} \cap  \overline{W(r_2)} =
W_0$ and $\overline{W(r_3)} \cup \overline{W(r_4)} = X_r, \mbox{   }$ 
$\overline{W(r_3)} \cap  \overline{W(r_4)} = W_1$.  We find that

{\footnotesize
\begin{eqnarray}\label{MV2}
&& 0 \rightarrow H_2(W_0) \xrightarrow{i_{1*} \oplus i_{2*}} H_2\left(\overline{W(r_1)}\right) \oplus H_2\left(\overline{W(r_2)}\right) \cr && \xrightarrow{j_{1*} - j_{2*}}  H_2(X_l^\infty)
\xrightarrow{\partial} H_1(W_0) \xrightarrow{i_{1*} \oplus i_{2*}} H_1\left(\overline{W(r_1)}\right) \oplus H_1\left(\overline{W(r_2)}\right) \rightarrow 0
\end{eqnarray}
}

\noindent
is exact,
where $i_1$ and $i_2$ are the inclusions of $W_0$ into $\overline{W(r_1)}$ and
$\overline{W(r_2)}$ and $j_1$ and $j_2$ are the inclusions of
$\overline{W(r_1)}$ and $\overline{W(r_2)}$ into $X_l^\infty$.  Slightly more work shows that
$\partial [V] = [S_0]$, where $[V]$ is the fundamental class of a vertical line in $X_l^\infty$ and
$[S_0]$ is the class of the invariant circle.

\vspace{.05in}

We have $H_2(X_r) = \mathbb{Z}^{\{[\mathbb{P}]\}}$ with $\partial([\mathbb{P}])
= [S_1]$.  Using that $H_i(X_r) = 0$ for $i\neq 2,0$,
we find that the map
\begin{eqnarray}\label{H2_ISOMORPHISM} H_2(W_1) &\xrightarrow{i_{3*} \oplus
i_{4*}}& H_2\left(\overline{W(r_3)}\right) \oplus
H_2\left(\overline{W(r_4)}\right) \end{eqnarray} \noindent is an isomorphism
and the sequence \begin{eqnarray}\label{MV4} 0 \rightarrow
\mathbb{Z}^{\{[\mathbb{P}]\}} \xrightarrow{\partial} H_1(W_1)
\xrightarrow{i_{3*} \oplus i_{4*}} H_1\left(\overline{W(r_3)}\right) \oplus
H_1\left(\overline{W(r_4)}\right) \rightarrow 0 \end{eqnarray} \noindent is
exact, where where $i_3$ and $i_4$ are the inclusions of $W_1$ into
$\overline{W(r_3)}$ and $\overline{W(r_4)}$.

\section{Morse Theory for $W_1$ and $W_0$}\label{MORSE}
In this section we prove that if there are parameter values $B$ for which $W_1
\cap C = \emptyset$, then $H_1(\overline{W(r_3)})$ and $H_1(\overline{W(r_4)})$
are trivial.  Our computer experiments suggest that such $B$ exist, but we have
not proven their existence.  This proves the second half of the second part of
Theorem \ref{MAIN_THM}, which we will finish proving in the following two
sections.

\vspace{.05in}

Consider the function $h:\mathbb{C} \times \mathbb{P} \rightarrow \mathbb{R}$
given by $h { x \choose y} = \left| \frac{x}{x-1} \right| $
which is chosen so that
\begin{eqnarray} \label{HEQN1}
h \left(N {x \choose y} \right) 
= \left| \frac{x^2}{x^2-2x+1} \right| =  h \left({x \choose y} \right) ^2.
\end{eqnarray}

We will consider the restriction of $h$ to the super-stable separatrices $W_0$
and $W_1$ and use it as a Morse function to study their topology.
Because $W_0$ and $W_1$ intersect the critical value parabola $C$ in
real-analytic sets, the following geometric description of the critical points
of $h$ makes sense:

\begin{prop}\label{MORSE_CRIT}
Let $K$ be the set of points in $C$ where $W_0 \cap C$ (or $W_1 \cap C$)
is parallel to the level curves of $h|C$.  Then
the set of critical points of $h$ on $W_0$ and $W_1$ is
$\bigcup_{k=1}^\infty N^{-k}(K)$.
\end{prop}

\noindent
{\bf Proof:}
Applying the chain rule to Equation \ref{HEQN1} we find:
\begin{eqnarray} \label{HEQN2}
Dh \left( N{x \choose y} \right) \cdot DN {x \choose y} = 2
 h{x \choose y}  \cdot Dh {x \choose y}.
\end{eqnarray}

Because $\left| h{x \choose y
} \right| = 0$ only when $x=0$, 
Equation \ref{HEQN2} gives that if $Dh {x \choose y} =0$
for a point ${x \choose y}$
on $W_0$ then either:
\begin{enumerate} 
\item {$Dh \left(
N{x \choose y} \right)=0$ giving that 
${x \choose y}$ is an inverse image
(possibly an $n$-th inverse image) of another critical point of $h$.  Or,}

\item{$DN {x \choose y}$ is singular and $Dh \left( N{x \choose y}\right)$ is
killed by $DN {x \choose y}$.}

\end{enumerate}
The condition in the second case says that $(x,y)$ is on the critical points
locus of $N$ and that the curve $W_0 \cap C$ is tangent to the level curves of
$h|C$ at $N(x,y)$.
\Endproof

\begin{figure}
\begin{center}
\scalebox{1.0}{
\begin{picture}(0,0)%
\epsfig{file=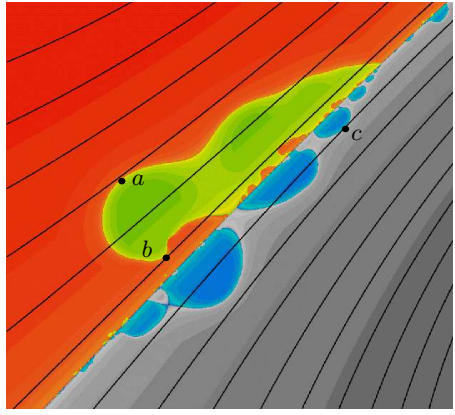}%
\end{picture}%
\setlength{\unitlength}{3947sp}%
\begingroup\makeatletter\ifx\SetFigFont\undefined%
\gdef\SetFigFont#1#2#3#4#5{%
  \reset@font\fontsize{#1}{#2pt}%
  \fontfamily{#3}\fontseries{#4}\fontshape{#5}%
  \selectfont}%
\fi\endgroup%
\begin{picture}(2800,2550)(1276,-2161)
\put(2063,-758){\makebox(0,0)[lb]{\smash{{\SetFigFont{9}{10.8}{\familydefault}{\mddefault}{\updefault}{\color[rgb]{0,0,0}$a$}%
}}}}
\put(3440,-476){\makebox(0,0)[lb]{\smash{{\SetFigFont{9}{10.8}{\familydefault}{\mddefault}{\updefault}{\color[rgb]{0,0,0}$c$}%
}}}}
\put(2127,-1212){\makebox(0,0)[lb]{\smash{{\SetFigFont{9}{10.8}{\familydefault}{\mddefault}{\updefault}{\color[rgb]{0,0,0}$b$}%
}}}}
\end{picture}%
}
\end{center}
\caption{Level curves of the Morse function $h$ in part of the critical value
parabola $C$. The points labeled $a$ and $b$ on $W_0$
are in $K$ and the point labeled $c$ on $W_1$
is in $K$.
Repeated inverse images of these points under $N$ are critical points of $h$
on $W_0$ and $W_1$.  Clearly
we have only
labeled a few of the points in $K$ that are visible.}
\label{FIG_MORSE}
\end{figure}

Notice that if $h:W_0 \rightarrow \mathbb{R}$ 
or if $h:W_1 \rightarrow \mathbb{R}$ has no critical points aside
from at $x=0$ or $1$, then the negative gradient flow $-\nabla h$ gives a deformation
retraction of $W_0$ to $S_0$ or the gradient flow $\nabla h$ gives a
deformation retraction of $W_1$ to $S_1$.

\begin{prop} If there are no points of intersection between $W_1$ and the
critical value parabola $C$, then $W_1$ is
homotopy equivalent to $S_1$.  
\end{prop}

\noindent {\bf Proof:} 
By Proposition \ref{MORSE_CRIT} if $W_1$ and $C$ are disjoint, there are no critical points of $h$.\Endproof

\begin{cor}
If there are no points of intersection between $W_1$ and the critical value parabola 
$C$, then the basins of attraction $W(r_3)$ and $W(r_4)$
for the roots $r_3 = (1,0)$ and $r_4 = (1,1-B)$ have trivial first 
and second homology groups. 
\end{cor}

\noindent {\bf Proof:}
This follows for the second homology from the isomorphism between $H_2(W_1)$
and $H_2(\overline{W(r_3)}) \oplus H_2(\overline{W(r_4)})$.  For the first
homology, $H_1(W_1) \cong \mathbb{Z}^{\{[S_1]\}} = {\rm Image}(\partial)$ and
exactness of  (\ref{MV4}) gives that
$H_1(\overline{W(r_3)})=0=H_1(\overline{W(r_4)})$.  Because $W_1$ is disjoint
from $C$, $B$ is not in the bifurcation locus, we can replace
$\overline{W(r_3)}$ and $\overline{W(r_4)}$ with the basins themselves.
\Endproof

\begin{prop}\label{W0_ALWAYSCRITICAL}
There are always critical points of $h:W_0 \rightarrow \mathbb{R}$.
\end{prop}

\noindent
{\bf Proof:}
Using implicit differentiation of $C(x,y)=0$, one can check that
there is a unique critical point of $h|C$ at the intersection
of $C$ with the line $Bx+2y-1=0$.  Since this line is the axis of symmetry
for $\tau$, it is in $W_0$.
\Endproof

Instead of studying the Morse function $h$ when $W_0$ or $W_1$ intersects $C$,
in the next two sections we will use linking numbers to show that such
intersections result in infinitely many loops corresponding to distinct generators of homology.

\section{Many loops in $W_0$ and $W_1$}\label{MANY_LOOPS}

Denote the vertical line in $X_l$ at a fixed value of $x$ by $\mathbb{P}_x$. 
Such vertical lines in $X_l$ correspond naturally to lines in
$X_l^\infty$ by means of the ``proper transform'' that is induced by the blow-up
operation.  

\vspace{.05in}
The Newton Map
$N$ maps $\mathbb{P}_x$ to $\mathbb{P}_{x^2/(2x-1)}$ by the rational map:
\begin{eqnarray*}
R_x(y) = \frac{y(Bx^2+2xy-Bx-y)}{(2x-1)(Bx+2y-1)}.
\end{eqnarray*}
\noindent
Notice that when $x=\frac{1}{B}$ and when $x=\frac{1}{2-B}$, a common term
cancels from the numerator and denominator of $R_x$, giving $R_x(y) =
\frac{y}{2}+\frac{1-B}{2(2-B)}$ and $R_x(y)=\frac{y}{2}$,  respectively.
The critical values of $R_x$ are the intersections of the critical value
parabola $C$ with the line $\mathbb{P}_{x^2/(2x-1)}$.  There are two distinct critical
values, except when $x=\frac{1}{B}$ or $\frac{1}{2-B}$.  

\begin{lem}
\label{NBHD} There are $\epsilon_0>0$ 
and $\epsilon_1>0$ so that if $|x-0| < \epsilon_0$, then $W_0 \cap 
\mathbb{P}_x$ forms a simple closed curve and so that if $|x-1| < 
\epsilon_1$, then $W_1 \cap \mathbb{P}_x$ forms a simple closed curve.
\end{lem}

\noindent
{\bf Proof:}
This is a direct consequence of the existence of $W_0^{loc}$ and $W_1^{loc}$ and
the fact that there is no possible recurrence in the dynamics for $x$ within
$X_l^\infty$ or $X_r$.
\Endproof

Most vertical lines $\mathbb{P}_x$ will be divided by $W_i$ ($i=0$ or $1$) into
exactly two simply connected domains.  However, if $W_i \cap C$ is non-empty in
any forward image of $\mathbb{P}_x$, then $W_i$ will divide $\mathbb{P}_x$ into many
more simply connected domains.  These are counted in the following proposition.

\begin{prop}\label{REGIONS} 
Let $\mathbb{P}_x$ be a vertical line whose $k$-th forward image
$\mathbb{P}_{\hat x}$ is divided by $W_i$ into exactly two simply connected
domains.  If  $W_i \cap C \neq \emptyset$ in $\mathbb{P}_{\hat x}$, then
$\mathbb{P}_x$ is divided by $W_i$ into between $2^k+2$ and $2^{k+1}$ simply
connected domains.  
\end{prop} 

We prove Proposition \ref{REGIONS} for $W_0$ in $X_l^\infty$ because the proof
is identical in $X_r$.  The following lemma and corollary are direct
consequences of the Riemann-Hurwitz Theorem.

\begin{lem} 
Let $R:\mathbb{P} \rightarrow \mathbb{P}$ be a ramified covering map of 
degree $d$ and let $U \subset \mathbb{P}$ be a simply connected open 
subset of $\mathbb{P}$ containing the image of at most one point of 
ramification of $R$.  Then, $R^{-1}(U)$ consists of a finite number of 
disjoint simply connected domains.
\end{lem}

The symmetry (\ref{SYMMETRY}) guarantees that there is at most one of the two
critical values of $R_{x}$ is in each simply connected domain, hence the
inverse image of each domain is a finite number of simply connected domains.
An easy check shows that if $U$ contains one of the critical values of $R_{x}$,
then $R_{x}^{-1}(U)$ is a single simply connected domain, while if $U$ does not
contain a critical value of $R_{x}$ it is two simply
connected domains.

\begin{cor}\label{ONESTEP}  
Let $\mathbb{P}_{x}$ be a vertical line whose image $\mathbb{P}_{\hat x}$ is
divided by $W_i$ into $m$ simply connected domains.  If $W_i \cap C
\neq \emptyset$ in $\mathbb{P}_{\hat x}$, then $W_i$ divides $\mathbb{P}_{x}$
into $2m$ simply connected domains. Otherwise it divides $\mathbb{P}_{x}$
into $2m-2$ simply connected domains.
\end{cor}

\noindent
The proof of Proposition \ref{REGIONS} follows from this Corollary and the fact
that there is a sufficiently high $k$ so that $|x_k| < \epsilon$ so that $W_0$ divides
$\mathbb{P}_{x_k}$ into exactly two domains.
\vspace{.05in}

\begin{figure}[!ht]
\begin{center}
\scalebox{.93}{
\begin{picture}(0,0)%
\epsfig{file=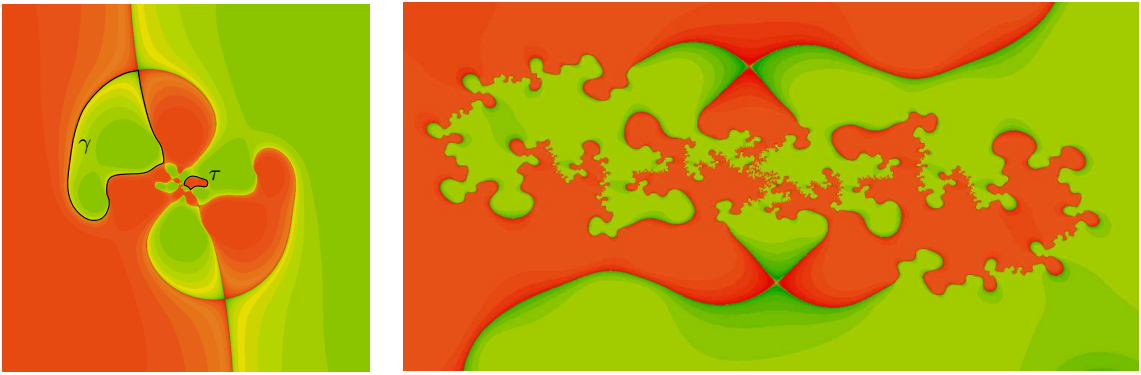}%
\end{picture}%
\setlength{\unitlength}{3947sp}%
\begingroup\makeatletter\ifx\SetFigFont\undefined%
\gdef\SetFigFont#1#2#3#4#5{%
  \reset@font\fontsize{#1}{#2pt}%
  \fontfamily{#3}\fontseries{#4}\fontshape{#5}%
  \selectfont}%
\fi\endgroup%
\begin{picture}(7674,2489)(665,-2468)
\put(1182,-968){\makebox(0,0)[lb]{\smash{{\SetFigFont{9}{10.8}{\familydefault}{\mddefault}{\updefault}{\color[rgb]{0,0,0}$\gamma$}%
}}}}
\put(2061,-1175){\makebox(0,0)[lb]{\smash{{\SetFigFont{9}{10.8}{\familydefault}{\mddefault}{\updefault}{\color[rgb]{0,0,0}$\tau$}%
}}}}
\end{picture}%

}
\end{center}
\caption{Left: a vertical line divided by $W_0$ into $10$ regions.  Loop
$\gamma$ has $\size(\gamma) = 1/10$ and loop $\tau$ has $\size(\tau) = 1/10$.
Right: a different vertical line that is divided by $W_0$ into $10$ regions.
This time, the loops bounding the regions are much more ornate.}
\end{figure}

Suppose that $W_0$ divides $\mathbb{P}_{x}$ into $2m$ simply connected domains
$U_1,\cdots,U_m$ in $W(r_1)$ and  $V_1,\cdots,V_m$ in $W(r_2)$.  Let $k$ be
chosen so that $W_0$ forms a simple closed curve in $\mathbb{P}_{\tilde{x}}$
(where $\tilde{x}$ is the $k$-th iterate of $x$ under $x \mapsto
\frac{x^2}{2x-1}$.) Denote by $U$ the domain in $\mathbb{P}_{\tilde{x}}$ within
$W(r_1)$ and by $V$ the domain in $\mathbb{P}_{\tilde{x}}$ within $W(r_2)$.

Under the mapping $N^k$, each domain $U_i$ covers $U$ with some degree $l_i$
and each domain $V_j$ covers $V$ with degree $p_j$.  Then: $\sum_{i=1}^m l_i =
2^k,  \qquad \sum_{i=1}^m p_i = 2^k$ because $U$ is covered by $\cup_{i=1}^m
U_i \subset \mathbb{P}_x$ with degree $2^k$.

For such a $U_i$ we will assign $\size(U_i) =
-\frac{l_i}{2^k}$ and such a $V_i$ we can assign 
$\size(U_i) = \frac{p_i}{2^k}$.
This is well defined because given $k_1$ and $k_2$ as above, the
$l_i$ corresponding to $k_1$ and the $l_i$ corresponding to $k_2$ will differ
by $2^{k_1-k_2}$.
\begin{eqnarray*}
\sum_{i=1}^m \size(U_i) = -1, \qquad \sum_{i=1}^m \size(V_i)=1
\end{eqnarray*}

In the next section we will see that $\size(U_i)$ for such a region equals the
linking number between $\gamma_i = \partial U_i$ and an appropriate geometric
object in $X_l^\infty$.

\section{Linking numbers}

\begin{figure}[!ht]\label{LINK}\begin{center}
\scalebox{.95}{
\begin{picture}(0,0)%
\epsfig{file=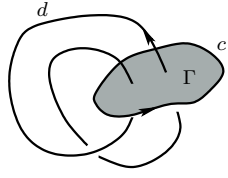}%
\end{picture}%
\setlength{\unitlength}{3947sp}%
\begingroup\makeatletter\ifx\SetFigFont\undefined%
\gdef\SetFigFont#1#2#3#4#5{%
  \reset@font\fontsize{#1}{#2pt}%
  \fontfamily{#3}\fontseries{#4}\fontshape{#5}%
  \selectfont}%
\fi\endgroup%
\begin{picture}(1577,1103)(1519,-966)
\put(2690,-393){\makebox(0,0)[lb]{\smash{{\SetFigFont{8}{9.6}{\familydefault}{\mddefault}{\updefault}{\color[rgb]{0,0,0}$\Gamma$}%
}}}}
\put(2912,-170){\makebox(0,0)[lb]{\smash{{\SetFigFont{8}{9.6}{\familydefault}{\mddefault}{\updefault}{\color[rgb]{0,0,0}$c$}%
}}}}
\put(1724, 53){\makebox(0,0)[lb]{\smash{{\SetFigFont{8}{9.6}{\familydefault}{\mddefault}{\updefault}{\color[rgb]{0,0,0}$d$}%
}}}}
\end{picture}%

}
\end{center}
\caption{Here  $lk(c,d) = +2$.}
\end{figure}

Classically one considers the linking number of two oriented loops $c$ and $d$
in $\mathbb{S}^3$.  The linking number $lk(c,d) \in \mathbb{Z}$ is found by
taking any oriented surface $\Gamma$ with oriented boundary $c$ and defining
$lk(c,d)$ to be the signed intersection number of $\Gamma$ with $d$ as in
Figure \ref{LINK}.  For this and many equivalent definitions of linking number
in $\mathbb{S}^3$ see \cite[pp. 132-133]{ROLF}, \cite[pp. 229-239]{BO_TU}, and \cite[Problems
13 and 14]{MILNOR_TOP}.

To see that this linking number is well-defined notice that assigning $lk(c,d) = [\Gamma]
\cdot [d]$, where $\cdot$ indicates the intersection product on
$H_*(\mathbb{S}^3,c)$, coincides with the classical definition.  If $\Gamma'$
is any other 2-chain with $\partial \Gamma ' = c$ then $\partial (\Gamma -
\Gamma') = [c]- [c] = 0$ and $(\Gamma - \Gamma')$ represents a homology class in
$H_2(\mathbb{S}^3)$.  Since $H_2(\mathbb{S}^3) = 0$, $[\Gamma - \Gamma'] = 0$
forcing $[\Gamma - \Gamma'] \cdot [d] = 0$.  Therefore: $[\Gamma] \cdot [d] =
[\Gamma'] \cdot [d]$, so that $lk(c,d)$ is well defined.

\vspace{.05in}
\noindent
{\bf Linking kernel: ${\cal L}Z_p(M)$}
\vspace{.05in}

Suppose that $M$ is a 3-dimensional manifold with $H_2(M) \neq 0$.  We can
define a linking number $lk(c,d)$ so long as the second argument $d$ has $[d]
\cdot \sigma = 0$ for every $\sigma \in H_2(M)$.  We define ${\cal L}Z_1(M)
\subset Z_1(M)$ to be the sub-module of one chains having this
property.  As before, given $d \in {\cal L}Z_1(M)$ the linking number $lk(c,d)$
is well-defined for $c$ disjoint from $d$ with $[c] = 0$.

Linking numbers work in a similar way if a manifold $M$ has dimension $m$:
one requires that $c$ and $d$ have dimensions summing to $m-1$ and one must
restrict to $c$ disjoint from $d$ with $[c] = 0$ and restrict to $d
\in {\cal L}Z_p(M)$, where $p$ is the dimension of $d$.

\subsection{Linking kernel for $X_l^\infty$}
Recall from Section \ref{COMP} that 
\begin{eqnarray*}
H_2(X_l^\infty) = \mathbb{Z}^{(L \cup \{[V]\})}
\end{eqnarray*}
\noindent
where $L$ is the set of exceptional divisors $E_i$ introduced in the sequence of
blow-ups and $V$ is the vertical line $x= \frac{1}{B(2-B)}$.

Recall from Proposition \ref{SELF_INTERSECTION}
that each exceptional divisor $[E_i]$ has
$[E_i]\cdot[E_i]\leq -1$ and that $[\vertline]\cdot[\vertline]=0$
so that if
$\omega = a_0 [\vertline] + a_1 [E_1] + \cdots + a_n [E_n]$
satisfies $\omega \cdot \sigma = 0$ for every $\sigma \in H_2(X_l^\infty)$ then
$a_i = 0$ for all $i \neq 0$. 

In summary, ${\cal L} Z_2(X_l^\infty) = \mathbb{Z}^{\{V\}}$.  The curves
$\gamma_i$ considered in the previous section each have linking number
$0$ with $\vertline$, since each $\gamma_i$ is entirely within some (other)
vertical line.  To show that any of the curves $\gamma_i$ are non-trivial
we will need to look for a different kind of object to link with.  We do this by
extending the definition of linking to closed currents.

\subsection{Generalities on currents}
Just as distributions are defined as the topological dual of smooth
functions with compact support, currents are the topological dual of
smooth differential forms with compact support.  

Let $A_c^{n-q}(M)$ denote the $(n-q)$-forms with compact
support on a smooth manifold $M$.  The linear maps  $T:A_c^{n-q}(M) \rightarrow
\mathbb{C}$ that are continuous are the {\bf currents of degree $q$} (or, as
some say, the currents of {\bf dimension $n-q$}) and are denoted by ${\cal
D}^q(M)$.  If $M$ has a complex structure, one defines the currents of
bi-degree $(p,q)$, denoted ${\cal D}^{p,q}(M)$ as the topological dual of the
$(n-p,n-q)$-forms with compact support $A_c^{n-p,n-q}(M).$
For more background on currents, consult \cite[section 3.1 and 3.2]{GH} or
the articles on complex dynamics \cite{HP2,SM1,SI1}.

An exterior derivative $d: {\cal D}^q(M) \rightarrow {\cal D}^{q+1}(M)$ is
defined as the adjoint to the classical exterior derivative on smooth forms
with compact support: $dT(\eta) = -T(d\eta)$.

On a complex manifold, one has two derivatives $\partial: {\cal D}^{p,q}(M)
\rightarrow {\cal D}^{p+1,q}(M)$ and $\overline{\partial}: {\cal D}^{p,q}(M)
\rightarrow {\cal D}^{p,q+1}(M)$, defined in the analogous way.  However, the
real operators $d = \partial + \overline{\partial}$ and $d^c =
\frac{i}{2\pi}(\overline{\partial}- \partial)$ are more often used in dynamics.
Currents that satisfy $dT = 0$ are referred to as {\em $d$-closed}.  We will
denote the $d$-closed currents of degree $q$ by $Z^q(M)$, and when desiring to
emphasize bi-degree, we will denote the $d$-closed currents of bi-degree $p,q$
by $Z^{p,q}(M)$.

Given a smooth form $\psi \in A^q(M)$, there is a current $T_\psi \in {\cal
D}^q(M)$ defined by $T_\psi(\eta) = \int_M \psi \wedge \eta$ for any $\eta
\in A_c^{n-q}(M)$.  Currents of this form are referred to as {\em smooth
currents}.  Using Stokes Theorem, one can check that $dT_\psi = T_{d\psi}$ so
that the inclusion $A_c^{n-q}(M) \rightarrow {\cal D}^q(M)$ given by $\psi
\mapsto T_\psi$ is a cochain map. 

A piecewise smooth, oriented $(n-q)$ chain $\Gamma$ in $M$ also defines a
current $T_\Gamma \in {\cal D}^q(M)$ given by $T_\Gamma(\eta) = \int_\Gamma
\eta$ for any $\eta \in A_c^{n-q}(M)$.  We will refer to currents that can be represented
this way as {\em currents of integration}.

We will often use work with {\bf closed, positive (1,1) currents}. These are
$(1,1)$ currents that are locally expressed as $T =  d d^c \phi$ for a
plurisubharmonic function $\phi$.  (See the $d d^c$-Poincar\'e Lemma.) We
denote the closed-positive $(1,1)$ currents on $M$ by $Z_+^{1,1}(M)$.

Typically one cannot pull back a current under a ramified mapping $F$.
One very special property of positive closed currents is that they can be pulled-back, using the potential
function: $F^*(\lambda) := dd^c (\phi \circ F)$, where $\phi$ is a
potential function for $\lambda$.  When $\lambda = T_\eta$ is a smooth current,
this pull-back coincides with the classical pull-back of smooth forms:
$F^*(T_\eta) =  T_{F^*(\eta)}$.


\subsection{Linking with currents}

The operator $d: {\cal D}^q(M) \rightarrow {\cal D}^{q+1}(M)$ satisfies $d
\circ d = 0$ and we denote the corresponding cohomology theory by
$H^*({\cal D}^{*}(M),d)$.  There is a natural map from the DeRham cohomology
$H^*_{DR}(M)$ into $H^*({\cal D}^{*}(M),d)$ induced by the inclusion of smooth forms
into the currents.

\begin{thm}{\bf (Approximation by smooth currents)}
The map $H^*_{DR}(M) \rightarrow H^*({\cal D}^{*}(M),d)$ is an isomorphism.  Furthermore,
the cohomology class of any closed current $L$ can be represented by a closed smooth form $\eta_L$
with support in an arbitrarily small neighborhood of the support of $L$.
\end{thm}

See \cite[pages 382-385]{GH} for a proof.  
\vspace{0.07in}

Given $T \in Z^{2}(M)$, and a piecewise smooth $2$-chain
$\sigma$ having $\partial \sigma$ disjoint from the support of $T$, 
there is a pairing:
\begin{eqnarray*}
C_2(M) \times Z^{2}(M) \rightarrow \mathbb{R}
\end{eqnarray*} 
\noindent
defined by $\left< \sigma,T \right> = \int_\sigma \eta_T$ 
were $\eta_T$ is a smooth form within the same cohomology class as $T$ with support is bounded away 
from $\partial \sigma$.  The existence of $\eta_T$ is garunteed by the approximation by smooth currents,
and the pairing is well defined since the integral depends only on the cohomology class of $\eta_T$.

When $T$ is a current of integration
integration over a piecewise smooth chain this pairing coincides with the usual
intersection number of piecewise smooth chains and when $T$ is given by a
smooth form, it coincides, by definition, with the standard pairing
$\int_\sigma T$.  (In fact, our pairing is a special case of the general intersection
number for closed currents of complimentary degrees \cite[p. 392]{GH} and \cite{DEMAILLY}.)

\begin{prop}\label{PAIRING_INV}
If $\lambda$ is a positive closed current and $F$ is a ramified mapping, we have
$\left< F_* \sigma, \lambda \right> = \left< \sigma, F^*\lambda \right>$.
\end{prop}

\noindent
{\bf Proof:}
Let $\eta_\lambda$ be a smooth approximation of $\lambda$ in the same cohomology class.
Then, $\left< F_* \sigma, \lambda \right> = \int_{F_* \sigma} \eta_\lambda = \int_\sigma F^* \eta_\lambda
= \left< \sigma, F^*\lambda \right>$, since $F^* \eta_\lambda$ is a smooth approximation of $F^*\lambda$.
\Endproof

\vspace{0.07in}

We define the linking kernel ${\cal L} Z^{2}(M)$ to be the space of closed
currents $T$ having $\left \langle \sigma, T \right \rangle = 0$ for every
$\sigma \in H_2(M)$.  Given $T \in {\cal L} Z^{2}(M)$, let $B_1^T(M)$ be the
$1$-boundaries in $M$ that are disjoint from the support of $T$.  We can define
a linking number with respect to $T$
\begin{eqnarray*}
lk(\cdot,T) : B_1^T(M) \rightarrow \mathbb{R}
\end{eqnarray*}

\noindent
by $lk(c,T) = \left \langle \Gamma ,T  \right \rangle$, where $\Gamma$ is any 2-chain with $\partial
\Gamma = c$.  Since $T \in {\cal L}Z^{2}(M)$, we have that $\left \langle \Gamma, T \right \rangle
= \left \langle \Gamma' ,T \right \rangle$ for any other $\Gamma'$ with $\partial \Gamma' = c$. 

\subsection{Finding an element of ${\cal L}Z^{2}(X_l^\infty)$}
In this subsection, we will find an element of ${\cal L}Z^{2}(X_l^\infty)$
by successively determining elements of ${\cal L}Z^{2}(X_l),$ ${\cal
L}Z^{2}(X_l^0),$ ${\cal L}Z^{2}(X_l^1),$ ${\cal
L}Z^{2}(X_l^2),\cdots$ where $X_k^j$ is the space $X_k$ after having
completed the blow-ups at level $j$.  In the limit, we will find an element of
${\cal L}Z^{2}(X_l^\infty)$, which in the next subsection will be useful
for linking.

Let $L_1$ be the invariant line that goes through $(0,0)$ and $(1,0)$, i.e.
$y=0$ and $L_2$ be the invariant line that goes through $(0,1)$ and $(1,1-B)$,
i.e. $y+Bx-1=0$.  (To remember the indexing, think that $L_1$ contains $r_1$
and $L_2$ contains $r_2$.) We can use the Poincar\'e-Lelong formula
(\cite[p. 388]{GH} or \cite{SM1}) to express the fundamental classes of these
lines as positive-closed currents: 
\begin{eqnarray*} [L_1] = \frac{1}{2\pi}
d d^c \log|y|,  \qquad [L_2] = \frac{1}{2\pi} 
d d^c \log|y+Bx-1|. \end{eqnarray*}

\noindent 
Both $L_1$ and $L_2$ intersect any given vertical line $\mathbb{P}$
with intersection number 1.  Because $[\vertline]$ is the sole generator of
$H_2(X_l)$ we have that $[L_2] - [L_1] \in {\cal L}Z^{2}(X_l)$.

Now, suppose that we want to find an element of ${\cal L}Z^{2}(X^0_l)$,
that is, a closed 2 current that evaluates to $0$ on every
element of $H_2(X^0_l) \cong \mathbb{Z}^{\{[\vertline],[E_p],[E_q]\}}$.  
In fact:
\begin{eqnarray*}
\left \langle E_p , [L_1] \right \rangle = 1 = \left \langle E_p, [L_2] \right \rangle  \mbox{     and     }
\left \langle E_q ,[L_1] \right \rangle = 0 = \left \langle E_q , [L_2] \right \rangle,
\end{eqnarray*}
\noindent
using standard intersection numbers for piecewise smooth chains,
so that $[L_2] - [L_1] \in {\cal L}Z^{2}(X^0_l)$.

This luck will not continue.  Let $z$ be one of the two 
preimages of $p$ that is in the invariant line $L_1$.  Since $L_1$
and $L_2$ intersect at the single point $p$, this forces that $z \notin L_2$.
Consequently:
$\left \langle E_z, [L_1] \right \rangle = 1 \neq 0 = \left \langle E_z, [L_2] \right \rangle$
so that $[L_2] - [L_1] \notin {\cal L}Z^{2}(X^1_l)$.

\vspace{0.1in}

We consider the $k$-th inverse images $N^{-k}(L_1)$
and $N^{-k}(L_2)$.
If we
denote by $N^k_1(x,y)$ and $N^k_2(x,y)$ the first
and second coordinates of $N^k$, then
the Poincar\'e-Lelong formula gives
\begin{eqnarray*}
[N^{-k}(L_1)] &=& \frac{1}{2 \pi} d d^c \log|N^k_2(x,y)|, \cr
[N^{-k}(L_2)] &=& \frac{1}{2 \pi} d d^c \log|N^k_1(x,y)+B\cdot N^k_2(x,y)-1|.
\end{eqnarray*}

\begin{lem}\label{EQUAL_ON_P} For every $k \geq 0$ we have 
\begin{eqnarray*}
\left \langle \vertline, [N^{-k}(L_1)] \right \rangle = \left \langle \vertline, [N^{-k}(L_2)] \right \rangle
\end{eqnarray*}
\end{lem}

\noindent
{\bf Proof:}
The $k$-th inverse images $N^{-k}(L_1)$ and $N^{-k}(L_2)$ both have
degree $2^k$ in $y$, so they each intersect a generic vertical line
transversely exactly $2^k$ times.  These intersection numbers coincide with the 
pairings.
\Endproof

\begin{prop}\label{EQUAL_LEVEL_K}
 $[N^{-(k+1)}(L_2)] - [N^{-(k+1)}(L_1)] \in {\cal L}Z^{2}(X^{k}_l)$
\end{prop}

\noindent
{\bf Proof}
Let $E_z$ be any one of the exceptional divisors in $X_l^k$.
Using Proposition \ref{EXCEPIONAL_COVERS},
there is some $d$ and some $l \leq k+1$ so that $N^{l}$ 
maps $E_z$ to $\vertline$
by a ramified cover of degree $d$ (possibly with $d = 0$.)  Just as
in the discussion above:
\begin{eqnarray*}
\left \langle E_z , [N^{-(k+1)}(L_1)] \right \rangle = \left \langle N^{l}(E_z), [N^{-(k+1)+l}L_1] \right \rangle =
        d \left \langle \vertline, [N^{-(k+1)+l} L_1] \right \rangle  \\
\left \langle E_z, [N^{-(k+1)}(L_2)] \right \rangle = \left \langle N^{l}(E_z), [N^{-(k+1)+l}L_2] \right \rangle =
        d \left \langle \vertline, [N^{-(k+1)+l} L_2] \right \rangle  \\
\end{eqnarray*}
\noindent
Here we are using Proposition \ref{PAIRING_INV} to obtain the first equality in each equation.  (One must check that
the Poincar\'e-Lelong equation gives that $[N^{-(k+1)}(L_i)] = (N^l)^* [N^{-(k+1)+l}L_i]$ for $l \leq k+1$.)
Then, Lemma \ref{EQUAL_ON_P} gives that the two terms on the right hand side of each equation are equal.

Since $H_2(X^k_l)$ is generated by the fundamental classes of $\vertline$
and the fundamental classes of each of the exceptional divisors $E_z$
we conclude that 
$[N^{-(k+1)}(L_2)] - [N^{-(k+1)}(L_1)] \in {\cal L}Z^{2}(X^k_l)$.
\Endproof

\vspace{.05in}

Since $X_l^\infty = \varprojlim (X_l^k,\pi)$ and $[N^{-(k+1)}(L_2)] -
[N^{-(k+1)}(L_1)] \in {\cal L}Z^{2}(X^k_l)$ we expect that a limit 
as $k \rightarrow \infty$ of $[N^{-(k+1)}(L_2)] - [N^{-(k+1)}(L_1)]$ will be
an element of ${\cal L}Z^{2}(X^\infty_l)$.  For such a limit to converge we
must normalize $[N^{-(k+1)}(L_2)]$ and $[N^{-(k+1)}(L_1)]$.
Dividing by the degrees, we define:
\begin{eqnarray*}
\lambda_1^k &=& \frac{1}{2^k} [N^{-k}(L_1)] = \frac{1}{2 \pi} d d^c \frac{1}{2^k} \log|N^k_2(x,y)|, \\
\lambda_2^k &=& \frac{1}{2^k} [N^{-k}(L_2)] = \frac{1}{2\pi} d d^c \frac{1}{2^k} \log|N^k_1(x,y)+B\cdot N^k_2(x,y)-1|
\end{eqnarray*}
\noindent
Let
\begin{eqnarray*}
\lambda_1 &=& \lim_{k \rightarrow \infty} \lambda_1^k = \frac{1}{2 \pi} d d^c \lim_{k \rightarrow \infty} \frac{1}{2^k} \log|N^k_2(x,y)|, \\
\lambda_2 &=&  \lim_{k \rightarrow \infty} \lambda_2^k = \frac{1}{2 \pi} d d^c  \lim_{k \rightarrow \infty} \frac{1}{2^k} \log|N^k_1(x,y)+B\cdot N^k_2(x,y)-1|.
\end{eqnarray*}

We will first check that these limits exist and define positive-closed 1-1
currents, and then we will show that $\lambda_2 - \lambda_1 \in {\cal
L}Z^{2}(X^\infty_l)$.

\begin{prop}\label{HP_CURRENTS}  The limits
\begin{eqnarray*}
G_1(x,y) &=& \lim_{k \rightarrow \infty} \frac{1}{2^k} \log|N^k_2(x,y)| \\
G_2(x,y) &=&  \lim_{k \rightarrow \infty} \frac{1}{2^k} \log|N^k_1(x,y)+B\cdot N^k_2(x,y)-1|
\end{eqnarray*}
converge and are plurisubharmonic functions in the basins of attraction $W(r_1)$
and $W(r_2)$, respectively.  Hence, $\lambda_1 = \frac{1}{2\pi} d d^c G_1(x,y)$
and $\lambda_2 = \frac{1}{2\pi} d d^c G_2(x,y)$ are positive closed 1-1
currents on $X_l^\infty$: $\lambda_1,\lambda_2 \in
Z_+^{1,1}(X_l^\infty).$
\end{prop}

\noindent
{\bf Proof:}
To see that $G_1(x,y)$ and $G_2(x,y)$ are well-defined and plurisubharmonic,
we will show that $G_1(x,y)$ and $G_2(x,y)$ coincide with the potential
functions that were described by Hubbard and Papadopol
in \cite[p. 21]{HP} and \cite{HP2}.  We will do this for $G_1(x,y)$, and leave
necessary modifications for $G_2(x,y)$ to the reader.

Supposing that $(0,0)$ is a root, Hubbard and Papadopol \cite{HP} consider the
limit
\begin{eqnarray*}
G_{HP}(x,y) = \lim_{k \rightarrow \infty} \frac{1}{2^k} \log ||N^k(x,y)||
\end{eqnarray*}

\noindent which they show converges to a plurisubharmonic function on the
basin of $(0,0)$.  The reader should notice that $G_{HP}$ does not depend on
the choice of the norm $|| \cdot ||$ used in the definition because any
two different norms on a finite dimensional vector space are equivalent by a
finite multiplicative constant, which is eliminated by the multiplicative
factor of $\frac{1}{2^k}.$ Therefore, we can use the supremum norm.

We will show that $G_1 = G_{HP}$ on 
$W(r_1)$, to see that $G_1$ is plurisubharmonic.  

If $|N^k_2(x,y)| \geq |N^k_1(x,y)|$ for all $(x,y)$ as $k \rightarrow \infty$,
then the supremum norm coincides with $|N^k_2(x,y)|$ giving $G_1(x,y) =
G_{HP}(x,y)$.  This condition is equivalent to the condition:
\begin{eqnarray} \label{G_ONE_VAR}
\lim_{k \rightarrow \infty} \frac{1}{2^k}\log
\left|\frac{N^k_2(x,y)}{N^k_1(x,y)}\right| \geq 0.
\end{eqnarray}
which will now show is a consequence of
a standard result from the dynamics of one complex variable.

In \cite{HP}, the authors perform blow-ups at each of the four roots, and
observe that the Newton map $N$ induces rational functions of degree $2$ on
each of the exceptional divisors $E_{r_1}, E_{r_2}, E_{r_3},$ and $E_{r_4}$.
Let's  compute the rational function $s: E_{r_1} \rightarrow E_{r_1}$.
In the coordinate chart $m=\frac{y}{x}$, the extension to $E_{r_1}$
is obtained by:
\begin{eqnarray*}
s(m) = \lim_{x \rightarrow 0} \frac{mx(Bx^2+2mx^2-Bx-mx)}{x^2(Bx+2mx-1)} 
= m(B+m)
\end{eqnarray*}
since $x=0$ on $E_{r_1}$.

Since condition (\ref{G_ONE_VAR}) is a limit, it suffices to check it in an
arbitrarily small neighborhood of the origin.  In a small enough neighborhood,
we can replace $\frac{N^k_2(x,y)}{N^k_1(x,y)}$ with $s\left(\frac{y}{x}\right)$ obtaining
\begin{eqnarray} \label{G_ONE_VAR_2}
\lim_{k \rightarrow \infty} \frac{1}{2^k}\log
\left|\frac{N^k_2(x,y)}{N^k_1(x,y)}\right|  = \lim_{k \rightarrow \infty} \frac{1}{2^k}\log |s^k(m)| = G_s(m).
\end{eqnarray}

\noindent where $G_s(m)$ is the standard Green's function from
one variable complex dynamics associated to the polynomial $s(m)$.  This
last equality is actually a delicate but well-known result that was proved
by Brolin \cite{BROLIN}.  A more friendly proof is available
in \cite[Section 9]{SM1}.

Having the last equality, it is a standard result, for example see Milnor
\cite{MILNOR} pages 95 and 96, that $G_s(m) = 0$ on the filled Julia set $K(s)$
and that $G_s(m) > 0$ outside of $K(s)$. 

This justifies the replacement of the supremum norm from $G_{HP}$ by
$|N^k_2(x,y)|$, and hence gives that $G_1(x,y) = G_{HP}(x,y)$.  
\Endproof

\begin{cor}\label{POTENTIAL_EQN}
Let $s:E_{r_1} \rightarrow E_{r_1}$ be the polynomial induced by
the Newton map $N$ and let $G_s:E_{r_1} \rightarrow \mathbb{R}$
be its Green's function.
In a sufficiently small neighborhood of $r_1$,
\begin{eqnarray*} 
G_1(x,y) = G_s \left(\frac{y}{x}\right) - \log \left|\frac{1}{x}\right|.
\end{eqnarray*}
\end{cor}

\noindent
{\bf Proof:}
This comes directly from the algebra:
\begin{eqnarray*}
G_1(x,y) &=& \lim_{k \rightarrow \infty} \frac{1}{2^k} \log|N^k_2(x,y)| 
= \lim_{k \rightarrow \infty} \frac{1}{2^k} \left(\log \left|\frac{N^k_2(x,y)}
{N^k_1(x,y)}\right| + \log |{N^k_1(x,y)}|\right) \cr
&=& G_s\left(\frac{y}{x}\right) + \lim_{k \rightarrow \infty} \frac{1}{2^k}
\log |{N^k_1(x,y)}| 
=  G_s\left(\frac{y}{x}\right) + \log|x| = G_s \left(\frac{y}{x}\right) - \log \left|\frac{1}{x}\right|
\end{eqnarray*}
\noindent
because $\frac{N^k_2(x,y)} {N^k_1(x,y)} \approx s \left(\frac{y}{x}\right)$ and
$N^k_1(x,y) = \frac{x^2}{2x-1} \approx -x^2$ near
$r_1$.
\Endproof

\subsection{Nice properties of $\lambda_2$ and $\lambda_1$:}
In this subsection, we will prove some of the useful properties if  $\lambda_2$
and $\lambda_1$.  We will finish the subsection by showing that $\lambda_2
-\lambda_1 \in {\cal L} Z^{2}(X_l^\infty)$.

\begin{lem}{\rm ({\bf Normalization})}\label{INT_TO_1}
Suppose that $\mathbb{P}_x$ is a vertical line that is divided into
exactly two simply connected domains $U\subset W(r_1)$ and $V\subset W(r_2)$ by $W_0$.  Then:
\begin{eqnarray*}
\left \langle V, \lambda_2 \right \rangle = 1 = \left \langle U, \lambda_1 \right \rangle \mbox{ and }
\left \langle U, \lambda_2 \right \rangle = 0 = \left \langle V, \lambda_1 \right \rangle
\end{eqnarray*}
\end{lem}
{\bf Proof:}
Because $N_2^k(x,y)$ and $B N_1^k(x,y)+N_2^k(x,y)-1$ are of degree
$2^k$ in $y$,
both $\lambda_1^k$ and $\lambda_2^k$ are normalized to that 
$\left \langle V, \lambda_1^k \right \rangle = 1$ and $\left \langle V, \lambda_2^k \right \rangle = 1$.  
Since the potentials for $\lambda_1^k$ and $\lambda_2^k$ converge
to the potentials for $\lambda_1$ and $\lambda_2$, we have
\begin{eqnarray*}
\left \langle U, \lambda_1 \right \rangle = \left \langle U, \lim_{k \rightarrow \infty} \lambda_1^k \right \rangle
= \lim_{k \rightarrow \infty} \left \langle U, \lambda_1^k \right \rangle = \lim_{k \rightarrow \infty} 1 = 1.
\end{eqnarray*}
\noindent
and similarly for $\lambda_2$.  The proof that $\left \langle U, \lambda_2 \right \rangle = 0 = \left \langle V, \lambda_1 \right \rangle$ is identical.
\Endproof

\begin{cor}\label{NORMALIZED}
Suppose that $\mathbb{P}_x$ is  vertical line, then $\left \langle \mathbb{P}_x,
\lambda_2 \right \rangle = 1 =  \left \langle \mathbb{P}_x,
\lambda_1 \right \rangle.$ 
\end{cor}

\vspace{0.1in}

It follows directly from the definitions of $\lambda_1$ and $\lambda_2$ that
$N^*(\lambda_1) = 2\cdot \lambda_1$ and $N^*(\lambda_2) = 2\cdot \lambda_2$.
(For example $N^*(\lambda_1) = \frac{1}{2\pi} dd^c \lim_{k \rightarrow \infty}
\frac{1}{2^k}  \log |N_1^k \circ N| = \frac{1}{2\pi} dd^c \lim_{k \rightarrow
\infty} \frac{1}{2^k}  \log |N_1^{k+1}| = 2\cdot \lambda_1$.) In combination
with Proposition \ref{PAIRING_INV}, this gives:

\begin{lem}{\rm({\bf Invariance})} \label{INVARIANCE}
Suppose that $\Gamma$ is a piecewise smooth $2$-chain, then 
\begin{eqnarray*}
\left \langle N(\Gamma), \lambda_1 \right \rangle = 2 \cdot \left \langle \Gamma, \lambda_1 \right \rangle  \qquad
\left \langle N(\Gamma), \lambda_2 \right \rangle = 2 \cdot \left \langle \Gamma, \lambda_2 \right \rangle
\end{eqnarray*}
\end{lem}

\vspace{0.1in}

\begin{prop} {\rm({\bf Support disjoint from $W_0$})} \label{SUPPORT}
There is a neighborhood $\Theta$ of $W_0$ in $X_l^\infty$ which
is disjoint from the support of $\lambda_1$ and $\lambda_2$.
\end{prop}

\noindent
{\bf Proof:}
By construction, $\lambda_1$ has support in $\overline{W(r_1)}$ and $\lambda_2$
has support in  $\overline{W(r_2)}$.  We will find a neighborhood, which we
also call $\Theta$, of $W_0$ in $\overline{W(r_1)}$ that is disjoint from the
support of $\lambda_1$.  Clearly similar methods will work in
$\overline{W(r_2)}$ and the desired neighborhood is the union of the two.

Recall from Corollary \ref{POTENTIAL_EQN} that
$G_1(x,y) = G_s \left(\frac{y}{x}\right) - \log\left|\frac{1}{x}\right|$,
where $G_s$ is the Green's function associated to the polynomial
$s:E_{r_1}\rightarrow E_{r_1}$ induced by $N$ at $r_1$.  Recall that
$s(m) = m(B+m)$ in the coordinates $m=\frac{x}{y}$ on $E_{r_1}$, so that
$m=\infty$ is a superattracting fixed point.  (This is the standard situation
for a quadratic polynomial.)

It is a standard result from one-variable dynamics, for example see \cite{MILNOR} p. 96, that $G_s$ is
harmonic outside of the Julia set $J(s)$.  In particular, $G_s$ is harmonic in
a neighborhood of $\infty$ (not including $\infty$).  
A related standard result that $G_s$ has the singularity 
\begin{eqnarray*}
G(m)= \log|m| + O(1) \mbox{ as } m \rightarrow \infty
\end{eqnarray*}
\noindent 
We check that this singularity exactly cancels with $ - \log\left|\frac{1}{x}\right|$
coming from
$G_1(x,y) = G_s \left(\frac{y}{x}\right) - \log\left|\frac{1}{x}\right|$:
\begin{eqnarray*}
G_1(x,y) &=& \log\left|\frac{y}{x}\right|  - \log\left|\frac{1}{x}\right| + O(1)
\mbox{  as  } \left|\frac{y}{x}\right| \rightarrow \infty \cr
&=& \log|y| + O(1) \mbox{  as  } \left|\frac{y}{x}\right| \rightarrow \infty
\end{eqnarray*} 
Therefore, $G_1(x,mx)$ is harmonic on a
neighborhood $U$ of $m=\infty$, including the point $\infty$.  Choose $\theta > 0$ so
that if $|m|>\theta$, then $G_1(x,mx)$ is harmonic.

Let $\Theta_0 = \{(x,y) \in \overline{W(r_1)} \mbox{ such that }
|\frac{y}{x}| > \theta \}$.  This is the open cone of points in $W(r_1)$ with slope
to the origin greater than $\theta$.  Since the invariant circle $S_0$ is above
$m=\infty$, $\Theta_0$ is a neighborhood of $S_0$ (within $\overline{W(r_1)}$.)

By construction,
$\Theta =  \bigcup_{n=0}^\infty N^{-n} (\Theta_0)$
will be invariant under $N$ and open.  Because $\Theta_0$ is disjoint from the
support of $\lambda_1$, the invariance properties for $\lambda_1$ from Lemma
\ref{INVARIANCE} give that all of $\Theta$ must be disjoint from the support of
$\lambda_1$.

Finally, since $\Theta_0$ contains a neighborhood of $S_0$, and both $W_0$ and
$\Theta$ are invariant under $N$, $\Theta$ forms an open neighborhood of $W_0$.
\Endproof

In fact, using the smooth approximation theorem, one can also choose the smooth
approximations of $\lambda_1$ and $\lambda_2$ to have support bounded away from
$W_0$.

\begin{cor}\label{ZERO_ON_W0}
Given any piecewise smooth chain $\sigma \in W_0$, we have that
$\left \langle \sigma, \lambda_1 \right \rangle = 0$ and $\left \langle \sigma, \lambda_2 \right \rangle = 0$.
\end{cor}

\begin{prop} $\lambda_1-\lambda_2 \in {\cal L} Z^{2}(X_l^\infty)$
\end{prop}

\noindent
{\bf Proof:}
This proof will be even simpler than the proof of Proposition \ref{EQUAL_LEVEL_K}
because we directly use the invariance of $\lambda_1$ and $\lambda_2$
shown in Lemma \ref{INVARIANCE}.

By Corollary \ref{NORMALIZED}, we have $\left \langle \vertline, \lambda_1
\right \rangle = \left \langle \vertline, \lambda_2 \right \rangle$.  Any
exceptional divisor $E_z$ was created during the blow-ups at some level $k$,
and using Proposition \ref{EXCEPIONAL_COVERS} there is some $l$ so that
$N^{\circ (k+1)}$ maps $E_z$ to $\vertline=\mathbb{P}_{1/(B(2-B))}$ by a
ramified covering mapping of degree $l$, (possibly $l = 0$).  Then:
\begin{eqnarray*}
\left \langle E_z, \lambda_1 \right \rangle = \frac{l}{2^{k+1}} \left \langle \vertline, \lambda_1 \right \rangle = 
\frac{l}{2^{k+1}} \left \langle \vertline, \lambda_2 \right \rangle = \left \langle E_z, \lambda_2 \right \rangle.
\end{eqnarray*}
Hence $\left \langle E_z, \lambda_2-\lambda_1 \right \rangle = 0$ for any exceptional divisor $E_z$.

Since an element of $H_2(X_l^\infty)$ is a linear combination of the
fundamental class $[\vertline]$ with a finite number of fundamental classes of
exceptional divisors $E_z$, we have shown that $\lambda_2 - \lambda_1 \in {\cal
L} Z^{2}(X_l^\infty).$ \Endproof

\subsection{$H_1(W_0)$ is infinitely generated.}
From Section \ref{MANY_LOOPS} we have infinitely many cycles $\gamma_i$ in
$W_0$ of arbitrarily small ``size,'' and we now have $(\lambda_2-\lambda_1) \in
{\cal L} Z^{2}(X_l^\infty)$ with which we can try to link them.

Since $H_1(X_l^\infty) = 0$, every $1$-cycle in $X_l^\infty$ is a
$1$-boundary in $X_l^\infty$.  In particular, $Z_1(W_0)  \subset
B_1(X_l^\infty)$.  By Lemma \ref{SUPPORT}, the support of $\lambda_2-\lambda_1$
is disjoint from $W_0$, giving that $Z_1(W_0) \subset
B_1^{\lambda_2-\lambda_1}(X_l^\infty)$.  Hence, we can restrict
$lk(\cdot,\lambda_2-\lambda_1)$ to $1$-cycles in $W_0$: 
\begin{eqnarray*}
lk(\cdot,\lambda_2-\lambda_1): Z_1(W_0) \rightarrow \mathbb{R} 
\end{eqnarray*}

\begin{prop}
For every $\gamma \in Z_1(W_0)$, $lk(\gamma,\lambda_2-\lambda_1)$ depends only on
$[\gamma] \in H_1(W_0)$.  
\end{prop}

\noindent
{\bf Proof:}
Suppose that $\gamma_1-\gamma_2 = \partial \sigma$, with $\sigma \in C_2(W_0)$.
Then, Corollary \ref{ZERO_ON_W0} gives that $\left \langle \sigma,
\lambda_2-\lambda_1 \right \rangle = 0$, hence $lk(\gamma_1,
\lambda_2-\lambda_1) = lk(\gamma_1, \lambda_2-\lambda_1)$.
\Endproof

\begin{prop}\label{LINK_RATIONAL}
The image of $lk(\cdot,\lambda_2-\lambda_1): H_1(W_0) \rightarrow \mathbb{R}$ is
contained in the rationals $\mathbb{Q}$.
\end{prop}
\noindent
{\bf Proof:}
Recall from Section \ref{MV} that there is an $\epsilon > 0$ for which $W_0$
restricted to $|x| < \epsilon$ is homeomorphic to the product $S_0 \times
\mathbb{D}_\epsilon$.  Because any $\gamma \in Z_1(W_0)$ is compact, there
exists a sufficiently high iterate $N^k$ so that $N^k(\gamma)$ lies within $|x|
< \epsilon$.  Then $[N^k(\gamma)] = n\cdot[S_0]$ for some appropriate $n$.  Using the invariance
property, this gives $lk(\gamma, \lambda_2-\lambda_1) = \frac{n}{2^k} \in
\mathbb{Q}$.  \Endproof

From here on we will write $lk(\cdot,\lambda_2-\lambda_1): H_1(W_0) \rightarrow \mathbb{Q}$.

\begin{prop}\label{LINK_SIZE} Suppose that $\gamma_i$ is a curve in a vertical line bounded
by a simply connected domain $U_i$.  Then:
$lk(\gamma_i,\lambda_2-\lambda_1) = \size(U_i)$,
where $\size(U_i)$ was defined in Section \ref{MANY_LOOPS}.
\end{prop}

\noindent
{\bf Proof of Proposition \ref{LINK_SIZE}:}
Recall that $\size(U_i)$ is defined as $\pm \frac{l_i}{2^k}$ where $k$ is such
that $N^k$ maps to a vertical line $\mathbb{P}_x$ that is divided by $W_0$
into only two domains $U \subset W(r_1)$ and $V \subset W(r_2)$ and where $l_i$ is the degree of this mapping
to $U$ or $V$.  The
sign is $-$ if $U_i$ is mapped to $U$ and $+$ if $U_i$ is mapped to $V$.
Without loss in generality, suppose that $U_i$ is mapped to $U$, and hence
$\size(U_i) < 0$.  Using Lemma \ref{INVARIANCE} we have that:
\begin{eqnarray*}
\left \langle U_i, \lambda_2-\lambda_1 \right \rangle = \frac{1}{2^k} \left \langle N^{k}(U_i), \lambda_2-\lambda_1 \right \rangle
=  \frac{1}{2^k} \left \langle l_i U, -\lambda_1 \right \rangle  &=& -\frac{l_i}{2^k} \left \langle U_i, \lambda_1 \right \rangle
= -\frac{l_i}{2^k} = \size(U_i)
\end{eqnarray*}

\noindent
where we are using that $\left \langle U, \lambda_2 \right \rangle = 0$ and $\left \langle U, \lambda_1 \right \rangle = 1$.
\Endproof

\begin{cor}\label{SMALL_ELEMENTS}
The image of the homomorphism
$lk(\cdot,\lambda_2-\lambda_1):H_1(W_0) \rightarrow \mathbb{Q}$
contains elements of arbitrarily small, but non-zero, absolute value.
\end{cor}

\noindent
This gives us our desired result:
\begin{cor}
The homology group $H_1(W_0)$ is infinitely generated.
\end{cor}

\noindent
Notice that an additive subgroup of $\mathbb{Q}$ that is dense must be infinitely
generated, but a dense additive subgroup of $\mathbb{R}$ typically is not infinitely generated
because the generators can be incommensurable.

\vspace{.05in}

Recall the Mayer-Vietoris exact sequence (\ref{MV2}):
{\small
\begin{eqnarray*}
H_2\left(\overline{W(r_1)}\right) \oplus H_2\left(\overline{W(r_2)}\right)
\rightarrow  H_2(X_l^\infty)  \xrightarrow{\partial} H_1(W_0) \rightarrow
H_1\left(\overline{W(r_1)}\right) \oplus H_1\left(\overline{W(r_2)}\right)
\rightarrow 0
\end{eqnarray*}
}
\noindent
If ${\rm Image}(\partial) =0$, or even if we knew that $|\size(\partial(\sigma))|$
were bounded away from $0$ for every $\sigma \in H_2(X_l^\infty)$,
we would be able to conclude that $H_1\left(\overline{W(r_1)}\right)$ and $H_1\left(\overline{W(r_2)}\right)$ are 
infinitely generated.  However, this is not the case.

\begin{prop} \label{ARB_SMALL_BDRY}
There are $\sigma \in H_2(X_l^\infty)$ with
$|lk(\partial(\sigma),\lambda_2-\lambda_1)|>0$ arbitrarily small.
\end{prop}

\noindent
{ \bf Proof:}
For every $k$, there exists some exceptional divisor $E$ having $N^{k}:E
\rightarrow V$ an isomorphism.  
For generic parameter
values $B\in S$, any exceptional divisor at a $(k-1)$-st
inverse image of $p$ will have this property, since, for generic
$B$ there is a single exceptional divisor above each point that we have blown
up, and $N:E_z \rightarrow E_{N(z)}$ is always an isomorphism.

For the values of $B \notin S$, there may be many blow-ups done at each
$(k-1)$-st inverse image of $p$.   We take a detailed look at the sequence
of blow-ups from Section \ref{SECTION_BLOWUPS} that was used to create
$X_l^{k-1}$ from $X_l^{k-2}$.  One must check that for each exceptional divisor
$E^{i}_{N(z)}$ that occurs in the sequence of blow-ups at $N(z)$, there is an
exceptional divisor in the sequence of blow-ups at $z$ that maps isomorphically
to $E^{i}_{N(z)}$.  Therefore, for any $k$, one can find an exceptional divisor
$E$ so that $N^{k-1}:E \rightarrow E_p$ is an isomorphism.
Since $N:E_p \rightarrow V$ is always an isomorphism, $E$ is the desired
exceptional divisor.  

Because $N^{k}$ maps $E$ isomorphically to $V$, it maps
$\partial([E])$ to $\partial([V])$.   
The invariance property from Lemma
\ref{INVARIANCE} gives that 
\begin{eqnarray*}
lk(\partial([E]),\lambda_2-\lambda_1) =
\frac{1}{2^k} lk(\partial([V]),\lambda_2-\lambda_1)=\frac{1}{2^k}.  
\end{eqnarray*}
\Endproof Proposition \ref{ARB_SMALL_BDRY}.

\subsection{$H_1\left(\overline{W(r_1)}\right)$ and $H_1\left(\overline{W(r_2)}\right)$ are infinitely generated.}
The following idea will allow us to
show that $H_1\left(\overline{W(r_1)}\right)$ and $H_1\left(\overline{W(r_2)}\right)$ are infinitely generated,
despite the fact that $|lk(\partial(\sigma),\lambda_2-\lambda_1)|$ can
be arbitrarily small, but non-zero, for $\sigma \in H_2(X_l^\infty)$.

\vspace{0.05in}
\noindent
{\bf Even and odd parts of Homology:}

Recall from Proposition \ref{SYMMETRY} that $N$ has a symmetry of reflection $\tau$
about the line $Bx+2y-1=0$ which exchanges the basins of attraction.
This $\tau$ induces an involution $\tau_*$ on $H_*(X_l^\infty),
 H_*(W_0)$, and $H_*(W(r_1)) \oplus H_*(W(r_2))$.  Every homology class
$\sigma$ will have $\tau_*^2(\sigma) = \sigma$ and consequently
the eigenvalues of $\tau$ are $\pm 1$.

We say that a homology class {\em $\sigma$ is even} if it is in 
the eigenspace of $\tau_*$ corresponding to eigenvalue $+1$, and we 
say that {\em $\sigma$ is odd} if it is in
the eigenspace of $\tau_*$ corresponding to eigenvalue $-1$.

Because the Mayer-Vietoris exact sequence commutes naturally with 
induced maps, we have a decomposition of the sequence (\ref{MV2})
into even and odd parts.  We will only need the odd part:
{\small
\begin{eqnarray*}
(H_2(\overline{W(r_1)}) \oplus H_2(\overline{W(r_2)}))^{\odd} \rightarrow 
H_2^{\odd}(X_l^\infty) \xrightarrow{\partial} H_1^{\odd}(W_0) \rightarrow
(H_1(\overline{W(r_1)}) \oplus H_1(\overline{W(r_2)}))^{\odd} \rightarrow 0
\end{eqnarray*}
}

\begin{lem} \label{TAU_ON_CURRENTS}
If $\sigma$ is some piecewise smooth chain, then: $\left \langle \sigma, \lambda_2 \right \rangle =
\left \langle \tau(\sigma), \lambda_1 \right \rangle$ and $\left \langle \sigma, \lambda_1 \right \rangle = \left \langle \tau(\sigma),
\lambda_2 \right \rangle$.
\end{lem}

\noindent
{\bf Proof:}

\noindent
Recall the definition of $\lambda_2$ and $\lambda_1$:
\begin{eqnarray*}
\lambda_1 &=& \frac{1}{2\pi} d d^c \lim_{k \rightarrow \infty} \frac{1}{2^k} \log|N^k_2(x,y)|, \\
\lambda_2 &=&  \frac{1}{2\pi} d d^c  \lim_{k \rightarrow \infty} \frac{1}{2^k} \log|N^k_1(x,y)+B\cdot N^k_2(x,y)-1|.
\end{eqnarray*}
\noindent
Since precomposition with $\tau$ exchanges the line $Bx+y-1=0$ with the line
$y=0$, Equation \ref{TAU_ON_CURRENTS} holds.  \Endproof

\begin{cor}
For every $[\gamma] \in H_1(W_0)$ we have:
$lk(\gamma,\lambda_2-\lambda_1) = -lk(\tau(\gamma),\lambda_2-\lambda_1)$.
\end{cor}

\noindent
{\bf Proof:}

\noindent
Suppose that $\sigma$ is a piecewise smooth $2$-chain with $\partial \sigma =
\gamma$.  Then we certainly have $\partial (\tau(\sigma)) = \tau(\gamma)$.
Lemma \ref{TAU_ON_CURRENTS} gives:
\begin{eqnarray*}
 lk(\gamma,\lambda_2-\lambda_1) = \left \langle \sigma, \lambda_2-\lambda_1 \right \rangle =
\left \langle \tau(\sigma), \lambda_1-\lambda_2 \right \rangle
= - \left \langle \tau(\sigma), \lambda_2-\lambda_1 \right \rangle =-lk(\tau(\gamma),\lambda_2-\lambda_1)
\end{eqnarray*}
\Endproof

\begin{prop}\label{ODD_BOUNDARY_ZERO}
If $\gamma \in H_1^{\odd}(W_0)$ is in the image of the boundary map 
$\partial:H_2^{\odd}(X_l^\infty) \rightarrow H_1^{\odd}(W_0)$,
then $lk(\gamma,\lambda_2-\lambda_1) = 0$. 
\end{prop}

\noindent
We first need the following lemma:
\begin{lem}\label{ANTI_COMMUTE}
For any exceptional divisor $E_z$ we have
\begin{eqnarray}\label{EQN_ANTI_COMMUTE}
\partial(\tau_*([E_z])) = - \tau_* (\partial([E_z]))
\end{eqnarray}
\end{lem}

\noindent
{\bf Proof:}
This proof will depend {\em essentially} on the explicit interpretation of the
boundary map $\partial$ from the Mayer-Vietoris sequence.  
In the following paragraph we closely paraphrase Hatcher \cite{HATCHER}, p. 150:

\vspace{.05in}
\noindent
The boundary map $\partial:H_n(X) \rightarrow H_{n-1}(A\cap B)$ can be made
explicit.  A class $\alpha \in H_n(X)$ is represented by a cycle $z$.  By
appropriate subdivision, we can write $z$ as a sum $x+y$ of chains in $A$ and
$B$, respectively.  While it need not be true that $x$ and $y$ are cycles
individually, we do have $\partial x = -\partial y$ since $z=x+y$ is a cycle.
The element $\partial \alpha$ is represented by the cycle $\partial x=-\partial
y$.
\vspace{.05in}

\begin{figure}[!ht]
\begin{center}
\begin{picture}(0,0)%
\includegraphics{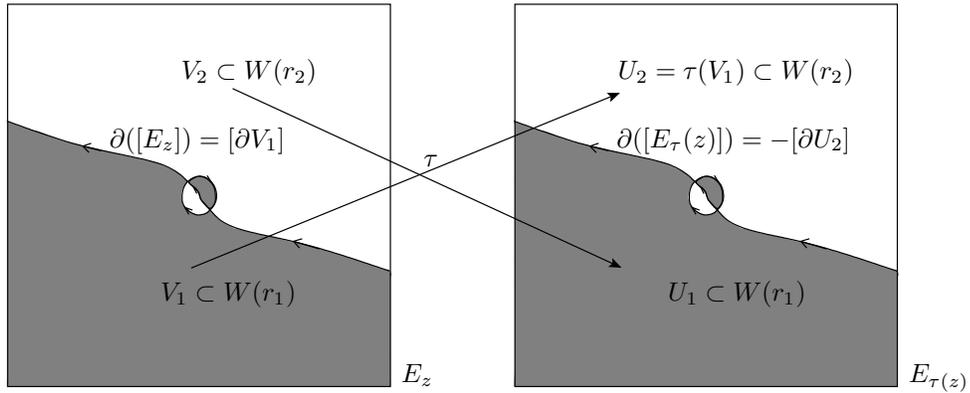}%
\end{picture}%
\setlength{\unitlength}{3947sp}%
\begingroup\makeatletter\ifx\SetFigFont\undefined%
\gdef\SetFigFont#1#2#3#4#5{%
  \reset@font\fontsize{#1}{#2pt}%
  \fontfamily{#3}\fontseries{#4}\fontshape{#5}%
  \selectfont}%
\fi\endgroup%
\begin{picture}(5676,2437)(1249,-2711)
\put(2226,-2149){\makebox(0,0)[lb]{\smash{\SetFigFont{10}{12.0}{\familydefault}{\mddefault}{\updefault}{\color[rgb]{0,0,0}$V_1 \subset W(r_1)$}%
}}}
\put(3740,-2664){\makebox(0,0)[lb]{\smash{\SetFigFont{10}{12.0}{\familydefault}{\mddefault}{\updefault}{\color[rgb]{0,0,0}$E_z$}%
}}}
\put(1906,-1184){\makebox(0,0)[lb]{\smash{\SetFigFont{10}{12.0}{\familydefault}{\mddefault}{\updefault}{\color[rgb]{0,0,0}$\partial([E_z])=[\partial V_1]$}%
}}}
\put(6925,-2664){\makebox(0,0)[lb]{\smash{\SetFigFont{10}{12.0}{\familydefault}{\mddefault}{\updefault}{\color[rgb]{0,0,0}$E_{\tau(z)}$}%
}}}
\put(5412,-2149){\makebox(0,0)[lb]{\smash{\SetFigFont{10}{12.0}{\familydefault}{\mddefault}{\updefault}{\color[rgb]{0,0,0}$U_1 \subset W(r_1)$}%
}}}
\put(5091,-1184){\makebox(0,0)[lb]{\smash{\SetFigFont{10}{12.0}{\familydefault}{\mddefault}{\updefault}{\color[rgb]{0,0,0}$\partial([E_\tau(z)])=-[\partial U_2]$}%
}}}
\put(3877,-1303){\makebox(0,0)[lb]{\smash{\SetFigFont{10}{12.0}{\familydefault}{\mddefault}{\updefault}{\color[rgb]{0,0,0}$\tau$}%
}}}
\put(5098,-753){\makebox(0,0)[lb]{\smash{\SetFigFont{10}{12.0}{\familydefault}{\mddefault}{\updefault}{\color[rgb]{0,0,0}$U_2 = \tau(V_1) \subset W(r_2)$}%
}}}
\put(2356,-753){\makebox(0,0)[lb]{\smash{\SetFigFont{10}{12.0}{\familydefault}{\mddefault}{\updefault}{\color[rgb]{0,0,0}$V_2 \subset W(r_2)$}%
}}}
\end{picture}

\end{center}
\caption{Showing that  $\partial(\tau_*[E_z]) = - \tau_* (\partial([E_z]))$.
}
\label{FIG_BOUNDARY}
\end{figure}

We use this explicit interpretation of $\partial$ to check Equation
\ref{EQN_ANTI_COMMUTE}.  Notice that $\tau_*{([E_z])}= [E_{\tau(z)}]$
consistent with the orientation 
that $E_z$ and $E_{\tau(z)}$ have as Riemann
surfaces.  Therefore  we have that $\partial(\tau_*([E_z]))= \partial([E_{\tau(z)}])=
[\partial U_1] = -[\partial
U_2]$, where $U_1$ is the oriented region of $E_{\tau(z)}$ that is in $\overline{W(r_1)}$
and $U_2$ is the oriented region of  $E_{\tau(z)}$ that is in $\overline{W(r_2)}$.

Similarly $\partial([E_z])=[\partial V_1] = -[\partial V_2]$, where $V_1$ and
$V_2$ are $E_z \cap \overline{W(r_1)}$ and $E_z \cap \overline{W(r_2)}$.  Because $\tau$ maps $E_z$
to $E_{\tau(z)}$  swapping $\overline{W(r_1)}$ with $\overline{W(r_2)}$ we have:
\begin{eqnarray*}
\tau_*(\partial([E_z])) = [\partial U_2] = -\partial(\tau_*([E_z]))
\end{eqnarray*}
\Endproof

\noindent
{\bf Proof of Proposition \ref{ODD_BOUNDARY_ZERO}:}

\noindent
Since elements of the form $[E_z]-[\tau(E_z)]$ span $H_2^{\odd}(X_l^\infty)$,
we need only check that the images of differences like this under $\partial$
have $0$ linking number:
\begin{eqnarray*}
 lk(\partial([E_z]-[\tau(E_z)]),\lambda_2-\lambda_1) &=&
lk(\partial([E_z])-\partial(\tau_*([E_z])),\lambda_2-\lambda_1) \cr &=&
lk(\partial([E_z])+\tau_*(\partial([E_z])),\lambda_2-\lambda_1) = 0
\end{eqnarray*}
\noindent
The last term is $0$ by Lemma \ref{ANTI_COMMUTE}.
\Endproof

\begin{prop}\label{ODD_SMALL_ELEMENTS}
The image of
$lk(\cdot,\lambda_2-\lambda_1): H_1^{\odd}(W_0) \rightarrow \mathbb{Q}$
contains elements of arbitrarily small, but non-zero absolute value.
\end{prop}

\noindent
{\bf Proof of Proposition \ref{ODD_SMALL_ELEMENTS}:}

\noindent
Recall from Proposition \ref{SMALL_ELEMENTS} that we can find $1$-cycles
$\gamma$ that have $lk(\gamma,\lambda_2-\lambda_1)$ arbitrarily small, but
non-zero.  Notice that $[\gamma - \tau(\gamma)]$ is obviously odd, and using
Lemma \ref{ANTI_COMMUTE}:
\begin{eqnarray*}
lk(\gamma-\tau(\gamma),\lambda_2-\lambda_1) &=& lk(\gamma,\lambda_2-\lambda_1)
- lk(\tau(\gamma),\lambda_2-\lambda_1) \cr &=&lk(\gamma,\lambda_2-\lambda_1) +
  lk(\gamma,\lambda_2-\lambda_1) = 2lk(\gamma,\lambda_2-\lambda_1).
\end{eqnarray*}
\noindent
\noindent
Hence, by choosing $\gamma$ so that $lk(\gamma,\lambda_2-\lambda_1)$ is
arbitrarily small, but non-zero, we can make
$lk(\gamma-\tau(\gamma),\lambda_2-\lambda_1)$ arbitrarily small, but non-zero
with $[\gamma-\tau(\gamma)] \in H_1^{\odd}(W_0)$.  \Endproof


Recall the last part of the exact sequence on the odd parts of homology:
\begin{eqnarray*}
\rightarrow H_2^{\odd}(X_l^\infty) \xrightarrow{\partial} H_1^{\odd}(W_0) \xrightarrow{i_{1*}\oplus i_{2*}}
\left(H_1\left(\overline{W(r_1)}\right) \oplus H_1\left(\overline{W(r_2)}\right)\right)^{\odd} \rightarrow 0
\end{eqnarray*}
\noindent
where $i_1$ and $i_2$ are the inclusions $W_0 \hookrightarrow
\overline{W(r_1)}$and $W_0 \hookrightarrow \overline{W(r_2)}$ respectively.

As a consequence of Proposition \ref{ODD_BOUNDARY_ZERO}, given any {\small $\eta \in
\left(H_1\left(\overline{W(r_1)}\right) \oplus H_1\left(\overline{W(r_2)}\right)\right)^{\odd}$} we can define
$lk(\eta,\lambda_2- \lambda_1) = lk(\gamma,\lambda_2- \lambda_1)$ for any
$\gamma \in H_1^{\odd}(W_0)$ whose image under $i_{1*}\oplus i_{2*}$ is $\eta$.
As a consequence of Proposition \ref{ODD_SMALL_ELEMENTS} we know that there are
{\small $\eta \in \left(H_1\left(\overline{W(r_1)}\right) \oplus H_1\left(\overline{W(r_2)}\right)\right)^{\odd}$} with
arbitrarily small $|lk(\eta,\lambda_2- \lambda_1)|$.  This proves the
desired result:

\begin{thm}
Let $\overline{W(r_1)}$ and $\overline{W(r_2)}$ be the closures in $X_l^\infty$ of the basins of attraction of
the roots $r_1=(0,0)$ and $r_2=(0,1)$ under the Newton Map $N$.  Then
$H_1\left(\overline{W(r_1)}\right)$ and $H_1\left(\overline{W(r_2)}\right)$ are infinitely generated. 
\end{thm}

\noindent
Recall also:
\begin{cor}
For parameter values $B \in \Omega_r$, we can replace
$\overline{W(r_1)}$ and $\overline{W(r_2)}$ with $W(r_1)$ and $W(r_2)$
finding that $H_1(W(r_1))$ and $H_1(W(r_2))$ are also infinitely generated.
\end{cor}

\subsection{Linking with currents in $X_r$}
Much of the work in the previous few subsections was to make linking numbers well
defined in $X_l^\infty$,  overcoming the indeterminacy from the fact that
$H_2(X_l^\infty)$ is infinitely generated.  Because $H_2(X_r) \cong
\mathbb{Z}^{\{[\mathbb{P}]\}}$ it is relatively easy to find elements in
${\cal L} Z_2(X_r)$.
However, we can just mimic the work from the previous sub-sections in an appropriate way.

The major difference is that in $X_l^\infty$ there is always an intersection of
$W_0$ with $C$ resulting in loops in $W_0$ of arbitrarily small size.  In
$X_r$, we must stipulate that an intersection of $W_1$ with $C$ exists before
proving that the homology is infinitely generated, because there appear to be
parameter values for which there is no intersection.

If we define $\lambda_3$ and $\lambda_4$ in a similar way as we defined $\lambda_1$ and
$\lambda_2$, then the following are proven in an easy way:

\begin{prop}
If $W_1$ intersects the critical value parabola $C$, then $H_1(W_1)$ is infinitely generated.
\end{prop}

Since there is only one generator of $H_2(X_r)$
this directly gives:

\begin{thm}
If $W_1$ intersects the critical value parabola $C$, then $H_1\left(\overline{W(r_3)}\right)$ and $H_1\left(\overline{W(r_4)}\right)$ are infinitely generated.
\end{thm}

\noindent
where $\overline{W(r_3)}$ and $\overline{W(r_4)}$ are the closures in $X_r$
of the basins of attraction of roots $r_3 = (1,0)$ and $r_4 = (1,1-B)$
under $N$.

\begin{cor}
For parameter values $B \in \Omega_{reg}$, we can replace
$\overline{W(r_1)}$ and $\overline{W(r_2)}$ with $W(r_1)$ and $W(r_2)$
finding that $H_1(W(r_1))$ and $H_1(W(r_2))$ are also infinitely generated.
\end{cor}

\vspace{.15in}
\noindent
This is the last part of the proof of Theorem \ref{MAIN_THM}.
$\Box$

\appendix

\section{Blow-ups of complex surfaces at a point.}\label{APP_BLOW_UPS}
Further material is available in \cite[pp. 182-189
and 473-478]{GH} and the introduction of \cite{HPV}.  

\vspace{.05in}
Suppose that $R : \mathbb{C}^2 \rightarrow \mathbb{C}^2$ has a point of indeterminacy
at $(0,0)$.  Blowing-up at $(0,0)$ produces a new space:
\begin{eqnarray}
\widetilde{\mathbb{C}}^2_{(0,0)} = \left\{(z,l) \in \mathbb{C}^2 \times \mathbb{P}^1 \mbox{ : } z \in l \right\}
\end{eqnarray}
\noindent
to which we can often find an extension $R: \widetilde{\mathbb{C}}^2_{(0,0)} \rightarrow
\mathbb{C}^2$ with no indeterminacy.  Here $\mathbb{P}^1$ is identified with
the space of directions through $(0,0)$ in $\mathbb{C}^2$.  There is a natural projection $\rho
: \widetilde{\mathbb{C}}^2_{(0,0)} \rightarrow \mathbb{C}^2$ given by
$\rho(z,l) = z$ and the set $E_{(0,0)} = \rho^{-1}((0,0))$ is referred to as the
{\em exceptional divisor.} A standard check shows that the blow-up is
independent of the choice of coordinates hence well defined on a complex surface
$M$ at a point $z$.

A rational map $R:\mathbb{C}^2 \rightarrow \mathbb{C}^2$ can be lifted to a new
rational mapping $\widetilde{R}:\widetilde{\mathbb{C}}^2_{(0,0)} - E_{(0,0)}
\rightarrow \mathbb{C}^2$ be defining $R(z,l) = R(z)$ for $z\neq 0$.  If the
indeterminacy in $R$ at $(0,0)$ was reasonably tame, $\widetilde{R}$ extends by
continuity to all of $E_{(0,0)}$.  Otherwise, there will be points of
indeterminacy of $\widetilde{R}$ on $E_{(0,0)}$ to which $\widetilde{R}$ cannot
be extended.  One can try further blow-ups at these points to resolve these new
points of indeterminacy.  The extension of $\widetilde{R}$ to $E_{(0,0)}$ is
analytic except at any new points of indeterminacy because $E_{(0,0)}$ is a
space of complex co-dimension 1.

\begin{prop} \label{BLOWUPHOM}
If $M$ is a complex surface and $z$ is any point in $M$, then the blow-up
$\widetilde{M}_z$ has 
$H_2(\widetilde{M}_z) \cong H_2(M) \oplus \mathbb{Z}^{(\{[E_z]\})}$ and
$H_i(\widetilde{M}_z) \cong H_i(M)$ for $i \neq 2$. 
\end{prop}

\noindent
The proof is an application of the Mayer-Vietoris sequence on homology and the fact
that $\widetilde{\mathbb{C}}^2_{(0,0)}$ has the homotopy type of $\mathbb{P}$.

Further analysis shows that the fundamental class $[E_z]$ of an exceptional
divisor has self-intersection number $-1$.  Meanwhile, blowing up a smooth
point on any complex curve $C$ decreases the self-intersection of its homology
class $[C]$ by one.  (See \cite{GH}, for proof.)

\section{Proof of Theorem \ref{BAIRE}}\label{AP_BAIRE}
Let $S \subset \Omega$ be the set of parameter values $B$ for which no inverse
image of the point of indeterminacy $p$ or the point of indeterminacy $q$ is in
the critical value locus $C$.  We are especially interested in $B \in S$
because the sequence of blow-ups from Section \ref{SECTION_BLOWUPS} is
especially easy to describe for these $B$.  

Theorem \ref{BAIRE} states that 
$S$ is generic in the sense of Baire's Theorem, i.e. uncountable and
dense in $\Omega$.  
The proof will follow as a corollary to:
\begin{Thm}({\bf Baire})
Let $X$ be either a complete metric space, or a locally compact Hausdorf space.
Then, the intersection of any countable family of dense open sets in $X$ is
dense.  
\end{Thm}
\noindent
See Bredon \cite{BREDON} for a proof of Baire's Theorem.

\vspace{.05in}
\noindent
{\bf Proof of Theorem \ref{BAIRE}:}
Let $S_n \subset \mathbb{C}$ be the subset of parameter values $B$ for which
none of the $n$-th inverse images of $p$ or $q$  under $N$ are in the critical
value locus $C$.

\begin{lem}\label{S_OPENDENSE}
 $S_n$ is a dense open set in $\mathbb{C}$ 
\end{lem}

\noindent
{\bf Proof:}
Let $R_n$ be the set of $B$ for which an $n$-th inverse image of $p$ is in $C$
and let $T_n$ be the set of $B$ for which an $n$-th inverse image of $q$ in
$C$.  We will show that $R_n$ and $T_n$ are finite, so that $S_n = \Omega
- (R_n \cup T_n)$ is a dense open set.

\begin{lem}\label{TN_FINITE}
$T_n$ is a finite set.
\end{lem}

\noindent
{\bf Proof:}
We have $B \in T_n$ if:
\begin{eqnarray} \label{CONDITION_T}
y^2+Bxy+\frac{B^2}{4}x^2-\frac{B^2}{4}x-y=0, \qquad
N^n_1(x,y) = \frac{1}{2-B}, \qquad
N^n_2(x,y) = \frac{1-B}{2-B}
\end{eqnarray}
\noindent
has a solution.  As always, $N^n_1$ and $N^n_2$
denote the first and second coordinates of $N^n$.  By clearing the denominators
in the second and third equations, condition \ref{CONDITION_T} can be expressed
as the common zeros of 3 polynomials $P_1(x,y,B)$, $P_2(x,y,B),$ and
$P_3(x,y,B)$ in the three variables $x, y,$ and $B$.  We will check there is no
common divisor of $P_1(x,y,B), P_2(x,y,B),$ and $P_3(x,y,B)$ so that the
solutions to \ref{CONDITION_T} form a finite set.

First, notice that $P_1(x,y,B) = y^2+Bxy+\frac{B^2}{4}x^2-\frac{B^2}{4}x-y$ is
irreducible.  It is sufficient to write an explicit biholomorphic map from
$\mathbb{C}^2$ to $\{P_1 =0 \} \subset \mathbb{C}^3$.  At a given $B$, the line
$Bx+2y=t$ intersects $\{P_1 =0 \}$ at a single point which we denote by
$f_B(t)$.  It is easy to check that $(t,B) \mapsto (f_B(t),B)$ provides the
desired isomorphism.

Hence $P_1$ has a factor in common with $P_2$ or $P_3$ if and only if $P_1$
divides $P_2$ or $P_3$.  We will show that this is impossible by examining
the lowest degree terms of $P_2$ and $P_3$.  If $P_1$ divides $P_2$ or
$P_3$, then the lowest degree term, $-y$, of $P_1$ must divide the lowest
degree term of $P_2$ or the lowest degree term of $P_3$.

We check by induction that the lowest degree term of $P_2$ is $\pm 1$ for
every $n$.  To simplify notation, let $a_k(x,y,B)$ be the polynomial
obtained by clearing the denominators from $N^k_1(x,y) = \frac{1}{2-B}$.

By clearing denominators of $N_1(x,y) = \frac{1}{2-B},$ we find $a_1(x,y,B) =
x^2(2-B) - 1(2x-1) = 2 x^2 -B x^2 -2x + 1$, so $a_1(x,y,B)$ has constant term
$\pm 1$.  Now suppose that $a_n(x,y,B)$ has constant term $\pm 1$.  By
definition, $a_{n+1}(x,y,B)$ is obtained by clearing the denominators of
$a_n(N_1(x,y),N_2(x,y),B) = 0$.  Because the denominators of both $N_1(x,y)$
and $N_2(x,y)$ have constant term $\pm 1$ and because $a_n(x,y,B)$ has
constant term $1$ we find that $a_{n+1}(x,y,B)$ has constant term $\pm 1$.

Because $P_2$ has constant term $\pm 1$ for every $n$ $P_1$ cannot divide
$P_2$, and we conclude that there are no common factors between $P_1$ and
$P_2$.

A nearly identical proof by induction shows that lowest degree term of $P_3$ is
also $\pm 1$ for each $n$.  Hence $P_1$ does not divide $P_3$, and we
conclude that $P_1$ and $P_3$ have no common divisors.

To see that $P_2$ and $P_3$ have no common divisors, notice that
$P_2(x,y,B) = 0$ is an equation for many disjoint vertical lines, while
$P_3(x,y,B) = 0$ stipulates that the $n$-th image of this locus has constant
$y=0$.  Since vertical lines are mapped to vertical lines by $N$, 
$P_2$ and $P_3$ can have no common factors.

Hence, $P_1, P_2,$ and $P_3$ are algebraically independent, so they have a
finite number of common zeros, giving that $T_n$ is a finite set.  \Endproof
Lemma \ref{TN_FINITE}.

\begin{lem}\label{RN_FINITE}
$R_n$ is a finite set.
\end{lem}

\noindent
{\bf Proof:}
Now we show that $R_n$, the set of $B$ so that an $n$-th inverse
image of $p$ under $N$ is in $C$, is finite.  In terms of equations, $R_n$ is the
set of $B$ so that:
\begin{eqnarray}
y^2+Bxy+\frac{B^2}{4}x^2-\frac{B^2}{4}x-y=0, \qquad
N^n_1(x,y) = \frac{1}{B}, \qquad
N^n_2(x,y) = 0 \label{CONDITION_R}
\end{eqnarray}
\noindent
has a solution.  Let $Q_1, Q_2,$ and $Q_3$ be the polynomials
equations resulting from clearing the denominators in Equation
\ref{CONDITION_R}.

The proof is the same as for $T_n$ except that a different proof is needed to
see that $Q_1$ does not divide $Q_3$.  An adaptation of the proof that $P_1$
does not divide $P_3$ fails because the lowest degree term of $Q_3$ has
positive degree in $y$.  We will check that $Q_1$ does not divide $Q_3$ and
leave the remainder of the proof to the reader.

The $x$-axis, $y=0$, is one of the invariant lines of $N$ and it intersects the
basins $W(r_1)$, $W(r_3)$ and the separator $\re(x) = 1/2$.  Therefore it is
disjoint from the two basins $W(r_2)$ and $W(r_4)$.  By definition,
$Q_3(x,y,B)$ is the equation for the $n$-the inverse image of the $x$-axis.
So, for a given $B$, the locus $Q_3(x,y,B)=0$ is also disjoint from the two
basins $W(r_2)$ and $W(r_4)$.

For every $B$, the critical value parabola $C$ goes through the four roots
$r_1$, $r_2$, $r_3$, and $r_4$, so it intersects all four basins of
attraction.  By definition,  $C$ is the zero locus $Q_1(x,y,B)=0$.  Therefore,
if $Q_1$ divides $Q_3$, there is a component of the zero locus $Q_3(x,y,B)=0$
intersecting all four basins $W(r_1)$, $W(r_2)$, $W(r_3)$ and $W(r_4)$ for
every $B$.
This is impossible, so $Q_1$ cannot divide $Q_3$.
\Endproof Lemma \ref{RN_FINITE} and
\Endproof Lemma \ref{S_OPENDENSE}.

Since $S_n$ is a dense open set in $\Omega$ for each $n$ and $S =
\cap_{n=0}^\infty S_n$, so it follows from Baire's Theorem that $S$ is
uncountable and dense in  the parameter space $\Omega$.  
\Endproof Theorem \ref{BAIRE}.

\vspace{.15in} \noindent {\bf \large Aknowledgements} 

This paper is a condensed version of my Ph.D. thesis for Cornell University.  I
thank the Department of Defense for generous financial support by means of a
National Defense Sciences and Engineering Graduate (NDSEG) fellowship and I
thank the National Science Foundation for further financial support by an
Interdisciplinary Graduate Education and Research Traineeship (IGERT)
fellowship.

My advisor John H. Hubbard suggested that I study the topology of the basins of
attraction for Newton's Method.  He provided mathematical guidance and
enthusiasm about this work and introduced to me the use of blow-ups to resolve
points of indeterminacy.  John Smillie and Eric Bedford provided a number of
suggestions including  encouraging me to learn about currents, which became a
key technique in this paper.  Allen Hatcher provided helpful comments about the
topology.  Alexey Glutsyuk provided very helpful discussions and suggestions
during the final stages of this project.

The computer program FractalAsm \cite{FA} written by  Karl Papadantonakis was
used to generate all of the images of the basins of attraction and was
invaluable for gaining an intuition about the topology of these basins.

The referee caught an error in my pairing between closed currents and piecewise
smooth chains (allowing for a correction)  and also provided many valuable
suggestions for streamlining and clarifying the final manuscript.  I thank the
referee greatly for his or her detailed consideration.

\bibliographystyle{plain}
\bibliography{newton.bib}

\end{document}